\documentclass[11pt]{article}
\usepackage{amsfonts}
\usepackage{color}
\usepackage{graphics}
\usepackage{cite}
\usepackage{latexsym}
\usepackage{amsmath}
\usepackage{amssymb}
\usepackage[dvips]{epsfig}
\usepackage{amscd}

\usepackage{color}
\usepackage{indentfirst}
\newtheorem{theorem}{Theorem}[section]
\newtheorem{proposition}{Proposition}[section]
\newtheorem{lemma}{Lemma}[section]

\newtheorem{remark}{Remark}[section]
\newtheorem{corollary}{Corollary}[section]

\def\pf{{\it Proof.}  }

\newcommand{\thatsall}{\hfill$\Box$}

\newcommand{\bt}{\begin{theorem}}
	\newcommand{\bl}{\begin{lemma}}
		\newcommand{\el}{\end{lemma}}
	\newcommand{\et}{\end{theorem}}

\newcommand{\bn}{\begin{eqnarray}}
	\newcommand{\en}{\end{eqnarray}}
\newcommand{\bnn}{\begin{eqnarray*}}
	\newcommand{\enn}{\end{eqnarray*}}

\newcommand{\bnnn}{\begin{eqnarray*}}
	\newcommand{\ennn}{\end{eqnarray*}}
\newcommand{\ben}{\begin{enumerate}}
	\newcommand{\een}{\end{enumerate}}

\newcommand{\ba}{\begin{aligned}}
	\newcommand{\ea}{\end{aligned}}
\newcommand{\be}{\begin{equation}}
	\newcommand{\ee}{\end{equation}}

\makeatletter
\@addtoreset{equation}{section}
\makeatother

\makeatletter
\@addtoreset{equation}{section}
\makeatother

\title{Global Strong Solution to Multi-Dimensional Compressible Quantum Navier-Stokes Equation in Bounded Domains with Large Initial Data}
\date{}
\author{Jiaxu Li$^{a}$, Boqiang L\"u$^{b}$, Huixiang Zhou$^{b}$, \thanks{email: jiaxvlee@gmail.com (J.X. Li),  lvbq86@163.com (B.Q. L\"u), zhouhx5010@163.com (H.X. Zhou) } \\
{\normalsize a. School of Mathematical Sciences,}\\
{\normalsize Shenzhen University, Shenzhen 518061, P. R. China;}\\
{\normalsize b. School of Mathematics and Computer Sciences}\\ {\normalsize \& Institute of Mathematics and Interdisciplinary Sciences,}\\ {\normalsize  Nanchang University, Nanchang 330031, P. R. China;}
}

\begin{document}
	\maketitle
	
	\begin{abstract}
		We consider the compressible quantum Navier-Stokes (QNS) system in a smooth bounded domain, supplemented with a curvature Navier-slip boundary condition for the velocity and a homogeneous Neumann boundary condition for the density. For arbitrarily large initial data away from vacuum, we establish the global existence and uniqueness of strong solutions when the adiabatic index
        \[ \gamma\in[1,\infty)\quad\text{for }N=2, \qquad \gamma\in\left[1,3\right)\quad\text{for }N=3. \] The argument relies on the effective velocity \(v=u+\nabla\log\rho\), which reveals the parabolic structure of the system. Compared with previous works in the periodic and Cauchy cases, additional difficulties arise from the treatment of boundary terms. We also use a bootstrap argument between the density integrability and the weighted integrability of the effective velocity to enlarge the admissible range of \(\gamma\) from $\left[1,11/4\right)$ to $\left[1,3\right)$.  Combining this bootstrap argument with a De Giorgi iteration for the density upper bound and a Nash-Moser type iteration for the positive lower bound, we derive a priori estimates on every finite time interval and extend the local strong solution globally in time. As far as we know, this is the first result on global strong solutions with large initial data for the initial-boundary value problem (IBVP) of the QNS system in bounded domains with  compatibility boundary conditions. 
	\end{abstract}

	


	\section{Introduction}

    The motion of compressible fluids with capillarity is described by the compressible Navier--Stokes--Korteweg (NSK) system:
\begin{equation}\label{x11}
		\left\{
		\begin{aligned}
			& \rho_t+\operatorname{div}(\rho u)=0,\\
			& (\rho u)_t+\operatorname{div}(\rho u\otimes u)
			+\nabla P(\rho)
			-\operatorname{div}\bigl(2\mu(\rho)D(u)\bigr)
			-\nabla\bigl(\lambda(\rho)\operatorname{div}u\bigr)
			=\operatorname{div}\mathbb K.
		\end{aligned}
		\right.
	\end{equation}
    Here, $t\ge 0$ denotes the time variable, $x\in\Omega\subset\mathbb R^N(N=2,3)$ is the spatial variable, \(\rho=\rho(x,t)\), \(u=u(x,t)\) and \(P(\rho)=\rho^\gamma\) are the density, velocity and pressure of the fluid,  respectively, and \(D(u)=\frac{\nabla u+(\nabla u)^{tr}}{2}\) is the deformation tensor. The viscosity coefficients \(\mu(\rho)\) and
	\(\lambda(\rho)\) are assumed to satisfy
	\[
	\mu(\rho)\ge 0,\qquad
	2\mu(\rho)+N\lambda(\rho)\ge 0.
	\]
	The Korteweg tensor \(\mathbb K\) is defined by
	\[
	\mathbb K=
	\left(
	\rho\operatorname{div}\bigl(k(\rho)\nabla\rho\bigr)
	-\frac{\rho k'(\rho)-k(\rho)}{2}|\nabla\rho|^2
	\right)\mathbb I
	-k(\rho)\nabla\rho\otimes\nabla\rho,
	\]
	where \(k(\rho)\) is the capillary coefficient and 
    \(\mathbb I\) denotes the identity matrix. 
	Throughout this paper, we impose the structural assumptions
	\begin{align}\label{x12}
		\mu(\rho)=\nu_0\rho,\qquad
		\lambda(\rho)=0,\qquad
		k(\rho)=\frac{k_0^2}{\rho},
	\end{align}
	where coefficients $\nu_0$ and $k_0$ satisfy
	\begin{align}\label{x13}
		\nu_0=k_0>0.
	\end{align}
	Under the assumptions \eqref{x12}--\eqref{x13}, system \eqref{x11} reduces to the compressible QNS system. By a normalization of the constants, we shall take
	\(\nu_0=k_0=1\).

    We study the IBVP of \eqref{x11} in a smooth bounded domain \(\Omega\subset\mathbb R^N\), \(N=2,3\), with the initial condition
	\begin{equation}\label{x14}
		\rho(x,0)=\rho_0(x),\qquad
		u(x,0)=u_0(x),\qquad x\in\Omega,
	\end{equation}
	and boundary condition
	\begin{equation}\label{x15}
		\partial_n\rho=0,\qquad
		u\cdot n=0,\qquad
		(D(u)n)_{tan}+\kappa(x)u_{tan}=0
		\quad\text{on }\partial\Omega.
	\end{equation}
    Where, $u_{tan}$ represents the projection of tangent plane of the vector $u$ on $\partial\Omega$, and $\kappa(x)$ is the curvature operator defined by
    \begin{equation}\label{x141}
    \kappa(x):T_x(\partial\Omega)\to T_x(\partial\Omega),
    \qquad
    \kappa(x)\xi
    :=
    D_\xi n(x)
    =
    (\xi\cdot\nabla)n(x).
    \end{equation} 
  
    To analyze the structure of the equation and make further estimates, we introduce the well-known effective velocity (see \cite{BreschDesjardinsLin,AntonelliSpirito2017,didier2019arma})
\begin{equation}\label{x161}
    v=u+\nabla\log\rho.
    \end{equation}
    Therefore, system \eqref{x11} can be rewritten as the following parabolic system 
    \begin{equation}\label{x16}
	\left\{
	\begin{aligned}
		&\rho_t+\operatorname{div}(\rho v)-\Delta\rho=0,\\
		& v_t+\bigl(v-2\nabla\log\rho\bigr)\cdot\nabla v
		-\Delta v+\nabla\Pi_\gamma(\rho)=0,
	\end{aligned}
	\right.
    \end{equation}
    where
    \[
    \Pi_\gamma(\rho)=
    \begin{cases}
	\frac{\gamma}{\gamma-1}\rho^{\gamma-1}, & \gamma>1,\\[2mm]
	\log\rho, & \gamma=1.
    \end{cases}
    \]
    Consequently, the initial effective velocity $v_0$ is defined by
    \begin{equation}\label{x17}
    v_0=u_0+\nabla\log\rho_0,
    \end{equation}
    Thanks to the boundary condition \eqref{x15}, the effective velocity $v$ satisfies the same Navier-type slip condition as $u$:
    \begin{equation}\label{x18}
	v\cdot n=0,\qquad
	(D(v)n)_{tan}+\kappa(x)v_{tan}=0
	\quad\text{on }\partial\Omega.
    \end{equation}
    The mathematical study of compressible fluids with capillarity dates back to Korteweg~\cite{Korteweg} and was further developed by Dunn-Serrin~\cite{DunnSerrin}. In the one-dimensional (1D) case, Germain-LeFloch~\cite{GermainLeFloch} established global finite-energy weak solutions and studied their convergence toward entropy solutions of the Euler system. Recently, Antonelli-Bresch-Spirito~\cite{AntonelliBreschSpirito2026} obtained global weak solutions for the 1D NSK system with power-law viscosity and capillarity coefficients, allowing vacuum without imposing an upper bound on the viscosity exponent. Global classical (strong) solutions with large initial data were also obtained under suitable assumptions on the viscosity and capillarity coefficients, see \cite{CharveHaspot,ChenChaiDongZhao2015,BurteaHaspotDegenerate} and the references therein. In multi-dimensions, Bresch-Desjardins-Lin~\cite{BreschDesjardinsLin} established global weak solutions allowing vacuum by deriving an additional entropy estimate. For strong solutions, however, the available results were mainly restricted to local existence or small initial data. Hattori-Li~\cite{HattoriLi1994,HattoriLi1996} established local strong solutions and global solutions for small non-vacuum initial data, while Danchin-Desjardins~\cite{DanchinDesjardins} obtained global smooth solutions for initial data close to a stable non-vacuum equilibrium. When it comes to the IBVP, Kotschote~\cite{Kotschote} further extended the local strong solution theory to domains with boundary.
    
    For the QNS system, which is a special case of the NSK system, J\"ungel~\cite{JungelEffective} established global smooth solutions with arbitrarily large initial data away from vacuum in 1D, which was later extended by Chen-Zhao~\cite{ChenZhaoQNS} to \(\nu>k\). In multi-dimensions, the global existence of weak solutions has been extensively studied, see \cite{JungelQNS,AntonelliSpirito2017,AntonelliSpirito2018,LacroixVioletVasseur,AntonelliHientzschSpirito} and the references therein. However, the global existence of strong solutions with arbitrarily large initial data remained a challenging problem, mainly due to the difficulty in obtaining a positive lower bound for the density. Very recently, there are some important progresses in \cite{HMZ,HuangGuLei2026,HuangMeng2026}. More precisely,  
    Huang-Meng-Zhang~\cite{HMZ} developed a modified Nash--Moser iteration based on the effective velocity, and established global strong solutions in periodic domains. Since their method cannot be applied directly to the Cauchy problem due to the nontrivial far-field behavior, Huang-Gu-Lei~\cite{HuangGuLei2026} used Littlewood--Paley theory and a refined De Giorgi iteration to obtain the corresponding whole-space result. For more general viscosity and capillarity coefficients, Gu-Huang-Meng-Zhou~\cite{GuHuangMengZhou2026} established global strong solutions with large initial data in periodic domains. Based on these results, Huang-Meng~\cite{HuangMeng2026} further obtained uniform-in-time upper and lower bounds for the density and studied the large-time behavior for the generalized system with $\alpha<1$. For the quantum case $\alpha=1$, they enlarged the three-dimensional (3D) range of $\gamma$ to $[1,11/4)$ and obtained a uniform-in-time upper bound for the density. However, the above results hold only for periodic or Cauchy problems. Hence, it is interesting to study the corresponding IBVP in smooth bounded domains.
    
    It should be mentioned here that similar problems have also been studied for other compressible fluid models. Inna-Saito~\cite{InnaSaito2023} established the local existence of strong solutions to the NSK system in general domains under Neumann and slip boundary conditions. Cai-Li~\cite{CaiLiBounded} and Cai-Li-L\"u~\cite{CaiLiLuExterior} respectively obtained global classical solutions to the compressible Navier--Stokes system in bounded and exterior domains with slip boundary conditions, provided that the initial energy is suitably small. However, these results do not provide the global estimates needed for the present QNS problem with arbitrarily large initial data. In particular, the density estimates based on the effective velocity in \cite{HMZ,HuangMeng2026} require suitable boundary conditions and additional treatment of the boundary terms in bounded domains.
    
    Therefore, this paper aims to establish a global strong solution theory for the IBVP in smooth bounded domains, while improving the admissible parameter range established in \cite{HuangGuLei2026} to make this theory applicable to more physical scenarios.
    
    
    
	
	Before stating the main results, we provide some notations that will be used throughout the entire paper. For any integrable function \(f\), we denote
	\[
	\int_{Q_T}f\,dxdt:=\int_0^T\int_{\Omega}f\,dxdt,
	\qquad
	Q_T:=\Omega\times[0,T].
	\]
	For \(1\le s\le\infty\), we use the standard notation
	\[
	L^s=L^s(\Omega),\qquad
	W^{k,s}=W^{k,s}(\Omega),\qquad
	H^s=W^{s,2}(\Omega).
	\]
	Throughout this paper, \(C\) denotes a generic positive constant that varies from line to line and depends on \(\Omega\), the initial data, and the parameters of the system. Its precise dependence will be specified when necessary.
	
	We now state our main result on the global existence of strong solutions to the IBVP.
	\begin{theorem}\label{th11}
    Let \(\Omega\subset\mathbb R^N\), \(N=2\) or \(N=3\), be a smooth bounded domain. Assume that the curvature operator \(\kappa(x)\) satisfies
	\begin{align}\label{x19}
		\kappa(x)\ge0,
		\qquad
		x\in\partial\Omega.
	\end{align}
	Assume that \(\gamma\) satisfies
	\begin{equation}\label{x111}
		\left\{\begin{aligned}
		&\gamma\in[1,\infty), \qquad \text{if } N=2\\
        &\gamma\in[1,3), \qquad \text{if } N=3
	\end{aligned}
		\right.
	\end{equation}
	and that the initial data \((\rho_0,u_0)\) satisfy
	\begin{align}\label{x113}
		0<\underline{\rho}_0\le \rho_0\le\overline{\rho}_0,
		\qquad
		\rho_0\in H^3(\Omega),
		\qquad
		u_0\in H^2(\Omega),
	\end{align}
	where \(\underline{\rho}_0\) and \(\overline{\rho}_0\) are positive constants. Moreover, assume that the initial data are compatible with the boundary conditions
	\eqref{x15}.
	Then the problem \eqref{x11}--\eqref{x15} admits a unique global
	strong solution \((\rho,u)\) satisfying for any \(0<T<\infty\) and \((x,t)\in\Omega\times[0,T]\),
	\begin{equation}\label{x115}
		\left\{
		\begin{aligned}
			& (C(T,\Omega))^{-1}\le \rho(x,t)\le C(T,\Omega),\\
			& \rho\in C([0,T];H^3(\Omega))\cap L^2(0,T;H^4(\Omega)),\\
			& \rho_t\in C([0,T];H^1(\Omega))\cap L^2(0,T;H^2(\Omega)),\\
			& u\in C([0,T];H^2(\Omega))\cap L^2(0,T;H^3(\Omega)),\\
			& u_t\in L^\infty(0,T;L^2(\Omega))\cap L^2(0,T;H^1(\Omega)).
		\end{aligned}
		\right.
	\end{equation}
	Here, the constant \(C(T,\Omega)>0\) depends on the initial data, \(T\) and \(\Omega\).
	\end{theorem}
    
    \begin{remark}\label{rm13}
	The boundary condition we introduced in \eqref{x15} is not only well-posed in mathematics, but also has certain physical significance. On the one hand, the homogeneous Neumann condition of density means that the density has no normal gradient at the wall, which is commonly associated with neutral wetting in diffuse-interface capillary models, see \cite{MagalettiGeorgoulasMarengo2020} for more information. On the other hand, the Navier-type velocity slip condition describes the stagnant layer of fluid near the wall, allowing fluid to slip when the shear stress is linearly related to the slip velocity, see \cite{CaiLiBounded,CaiLiLuExterior,Maxwell1879,NavierFluid1827} for more details.
    
    
    \end{remark}

    \begin{remark}\label{rmk12}
		When \(N=2\), the curvature operator \(\kappa(x)\) reduces to the
		scalar curvature determined by
		\[
		\kappa(x)\xi=\kappa\xi, \qquad \xi\in T_x(\partial\Omega).
		\]
		When \(N=3\), \(\kappa(x)\) is the shape operator of
		\(\partial\Omega\). Therefore, the condition \(\kappa(x)\ge0\) in
		\eqref{x19} is understood in the scalar sense for \(N=2\), while for
		\(N=3\) it means
		\[
		\xi\cdot\kappa(x)\xi\ge0,
		\qquad
		\forall\,\xi\in T_x(\partial\Omega).
		\]
	\end{remark}
    
    \begin{remark}\label{Rmk11}
		The boundary condition in \eqref{x15} for the velocity is equivalent to
		\begin{equation}\label{x151}
			u\cdot n=0,\qquad
			\operatorname{curl}u\times n=0
			\quad\text{on }\partial\Omega,
		\end{equation}
		by the definition of $\kappa(x)$ (see \cite[Remark~1.1]{CaiLiBounded} for details). For the same reason,
        the boundary condition in \eqref{x18} for the effective velocity is equivalent to
        \begin{equation}\label{x152}
        v\cdot n=0,\qquad \operatorname{curl}v\times n=0 \quad\text{on }\partial\Omega.
        \end{equation}
    \end{remark}
    
    \begin{remark}\label{rmk11}
		Compared with previous works \cite{HuangGuLei2026,HuangMeng2026}, we  enlarge the admissible range of \(\gamma\) in \eqref{x111}\(_2\) from $\left[1,11/4\right)$ to $\left[1,3\right)$ in 3D.
	\end{remark}

    We offer some comments on our analysis of this paper. The passage from the periodic and whole space settings of \cite{HuangGuLei2026,HuangMeng2026,HMZ} to a bounded domain presents fundamental challenges, which require the development of new technical approaches. 

    First, the choice of boundary conditions is crucial to our analysis. In order to make full use of the parabolic structure of \eqref{x16}, we need suitable boundary conditions for both $u$ and $v$. However, a general Navier-type slip boundary condition for $u$ (see \cite{NavierFluid1827,CaiLiBounded}) is not necessarily preserved under the transformation \eqref{x161}. To overcome this difficulty, we impose the boundary condition \eqref{x15}, whose velocity part is equivalent to \eqref{x151} (see \cite[Remark~1.1]{CaiLiBounded}). It follows from \eqref{x151} that $v$ satisfies \eqref{x152}, and hence the same Navier-type slip condition \eqref{x18}. Thus the boundary conditions are compatible with $v$, which allows us to work with \eqref{x16}.

    Next, the boundary terms arising from integration by parts cause additional difficulties in deriving the basic energy estimates. Thanks to \eqref{x15} and \eqref{x18}, the terms involving $u\cdot n$, $v\cdot n$, and $\partial_n\rho$ vanish. However, the remaining boundary integrals still need to be controlled. More precisely, using the boundary condition \eqref{x15}, the identity \eqref{x3121}, and the assumption \eqref{x19}, we have
	\[
	\begin{aligned}
		\mathcal B_{u}(t)
		&=2\int_{\partial\Omega}\rho\,\kappa(x)u_{tan}\cdot u_{tan}\,dS\ge0,\\
		\mathcal B_{v}(t)
		&=\int_{\partial\Omega}\rho\,\kappa(x)v_{tan}\cdot v_{tan}\,dS\ge0.
	\end{aligned}
	\]
	Additionally, the same argument shows that the boundary integral arising in the weighted estimate for $v$ is also nonnegative (see \eqref{x44} for details). Thus, we can obtain the basic energy estimates \eqref{x31} and \eqref{x39}, together with the weighted estimates for $v$, which are needed to establish the upper and positive lower bounds for the density.

    Finally, in order to extend the local strong solution globally in time, we need to derive the higher-order estimates. However, the elliptic estimates used in the periodic and Cauchy cases (see \cite{GuHuangMengZhou2026,HuangGuLei2026,HuangMeng2026,HMZ}) cannot be applied directly to the bounded domain problem, since the boundary conditions must be taken into account. To overcome this difficulty, we first use the Neumann elliptic estimate for $\rho$, which gives
	\begin{equation}\label{x156}
		\|\rho\|_{W^{k+2,p}(\Omega)}
		\le C\left(
		\|\Delta\rho\|_{W^{k,p}(\Omega)}
		+\|\rho\|_{L^p(\Omega)}
		\right),
	\end{equation}
	On the other hand, since $v$ satisfies the slip boundary condition \eqref{x152}, we obtain the following elliptic estimate:
	\begin{equation}\label{x157}
		\|v\|_{W^{k+2,p}(\Omega)}
		\le C\left(
		\|\Delta v\|_{W^{k,p}(\Omega)}
		+\|v\|_{L^p(\Omega)}
		\right),
	\end{equation}
	These estimates (see Lemma~\ref{lm27}--\ref{lm28} for details) play an important role in controlling the higher derivatives of $\rho$ and $v$. In addition, when estimating $v_t$, one needs to deal with the boundary integrals. Using the boundary identity \eqref{x3111}, the classical trace theorem, and Young's inequality, we obtain \eqref{x519}. Thus, we can establish the higher-order estimates in Section~5 and then extend the local strong solution globally in time.

    Moreover, the admissible range of $\gamma$ in \eqref{x111}\(_2\) is enlarged by a bootstrap argument (see Proposition~\ref{pr33} for details). Motivated by \cite{HuangMeng2026}, we first use weighted $L^p$-estimates for $v$ with $p>2$ to improve the integrability of the density, instead of directly obtaining the weighted $L^{q+2}$-estimate with $q>1$ as in \cite{HuangGuLei2026}. Indeed, starting from $Q_0=3$ (see Corollaries~\ref{cor31}--\ref{cor32}), we obtain the following iteration:
	\begin{equation}\label{x153}
		Q_{k+1}=Q_k+3-\gamma.
	\end{equation}
	Thus, for $2<\gamma<3$, the density integrability is improved at each step, which finally yields the weighted estimate needed for the De Giorgi iteration. Unfortunately, when $\gamma=3$, the gain vanishes and the bootstrap argument cannot proceed. This explains the restriction $\gamma<3$ in \eqref{x111}\(_2\).
	
    The rest of this paper is organized as follows. In Section 2, we collect the auxiliary results required for the bounded-domain analysis. Section 3 is devoted to the upper bound of the density, where the bootstrap argument is carried out in the bounded-domain setting. 
    Section 4 establishes the positive lower bound of the density. 
    Section 5 derives the higher-order estimates with the elliptic estimates. Finally, Section 6 completes the proof of Theorem~\ref{th11} by a continuation argument.
	
	\section{Preliminaries}
	
	In this section, we collect several elementary inequalities which will be used frequently in the following sections. The first is the space-time interpolation inequality, which can be found in \cite[Theorem~1.49]{NovotnyStraskraba}.
	
	\begin{lemma}\label{lm21}
		Let \(0<T<\infty\). Assume that \(0\le \theta_i\le 1\), \(i=1,\cdots,m\), and \( \sum_{i=1}^m \theta_i=1. \)
		Let \(1\le p_i,q_i\le\infty\), \(i=1,\cdots,m\), and define \(p,q\) by \[ \frac{1} {p}=\sum_{i=1}^m\frac{\theta_i}{p_i}, \qquad \frac{1}{q}=\sum_{i=1}^m\frac{\theta_i}{q_i}. \] Then \[ \|f\|_{L^p(0,T;L^q(\Omega))} \le \prod_{i=1}^m \|f\|_{L^{p_i}(0,T;L^{q_i}(\Omega))}^{\theta_i}. \]
	\end{lemma}

   According to \cite{EvansPDE}, in the follow, we state the Space-time H\"older inequality.
    
	\begin{lemma}\label{lm22}
		Let \(0<T<\infty\). Assume that \(1\le p_i,q_i\le\infty\), \(i=1,\cdots,m\), and define \(p,q\) by
		\[
		\frac{1}{p}=\sum_{i=1}^m\frac{1}{p_i},
		\qquad
		\frac{1}{q}=\sum_{i=1}^m\frac{1}{q_i}.
		\]
		Then
		\[
		\left\|\prod_{i=1}^m f_i\right\|_{L^p(0,T;L^q(\Omega))}
		\le
		\prod_{i=1}^m
		\|f_i\|_{L^{p_i}(0,T;L^{q_i}(\Omega))}.
		\]
	\end{lemma}

    We also recall the following Chebyshev inequality from \cite[Remark~1.18]{AmbrosioFuscoPallara}, which will be used in the De Giorgi iteration for the density upper bound.
    
	\begin{lemma}[Chebyshev inequality]\label{lm23}
		Let \(f\in L^p(\Omega)\) for some \(p>0\). Then for every \(\lambda>0\),
		\[
		|\{x\in\Omega:\ |f(x)|>\lambda\}|
		\le
		\frac{1}{\lambda^p}\int_\Omega |f(x)|^p\,dx.
		\]
	\end{lemma}

    Next, the well-known Gagliardo-Nirenberg’s inequality on bounded domains (see \cite[Theorem~1.2]{LiZhang2022}) will be used frequently later.
	\begin{lemma}[Gagliardo--Nirenberg inequality]\label{lm24}
		Let $\Omega\subset\mathbb R^N$, $N=2,3$, be a smooth bounded domain. For $1\le q<\infty$, $r>N$ and $q\le p<\infty$, there exists a constant $C>0$, depending only on $\Omega,N,p,q,r$, such that for any
		$f\in L^q(\Omega)\cap D^{1,r}(\Omega)$,
		\[
		\|f\|_{L^p(\Omega)}
		\le C\|f\|_{L^q(\Omega)}^{1-\theta}
		\|\nabla f\|_{L^r(\Omega)}^\theta
		+C\|f\|_{L^q(\Omega)},\qquad
		\theta=\frac{Nr(p-q)}{p\big(Nr+q(r-N)\big)}.
		\]
		Moreover, for any $g\in L^q(\Omega)\cap D^{1,r}(\Omega)$,
		\[
		\|g\|_{L^\infty(\overline\Omega)}
		\le C\|g\|_{L^q(\Omega)}^{\frac{q(r-N)}{Nr+q(r-N)}}
		\|\nabla g\|_{L^r(\Omega)}^{\frac{Nr}{Nr+q(r-N)}}
		+C\|g\|_{L^q(\Omega)}.
		\]
	\end{lemma}
	
	Finally, we state two elliptic estimates needed for the higher-order analysis. The first is a standard Neumann estimate; see \cite[Theorem~4.27]{NovotnyStraskraba}. It will be used in Sections~5 to control higher derivatives of the density.
	
	\begin{lemma}\label{lm27}
		Let \(1<p<\infty\), \(k=0,1,\cdots\), and let
		\(\Omega\subset\mathbb R^N\), \(N=2,3\), be a bounded smooth domain.
		Assume that \(\rho\) is sufficiently smooth and satisfies
		\[
		\partial_n\rho=0
		\quad\text{on }\partial\Omega.
		\]
		Then
		\begin{align}\label{x22}
			\|\rho\|_{W^{k+2,p}(\Omega)}
			\le
			C
			\left(
			\|\Delta \rho\|_{W^{k,p}(\Omega)}
			+
			\|\rho\|_{L^p(\Omega)}
			\right).
		\end{align}
		In particular, when \(p=2\),
		\begin{align}\label{x23}
			\|\rho\|_{H^{k+2}(\Omega)}
			\le
			C
			\left(
			\|\Delta \rho\|_{H^k(\Omega)}
			+
			\|\rho\|_{L^2(\Omega)}
			\right).
		\end{align}
	\end{lemma}

    In order to control higher derivatives of the effective velocity in Section~5, we need the Agmon-Douglis-Nirenberg (ADN) elliptic estimates \cite{ADN1964}.
	Consider the following elliptic system 
	\[ \begin{cases} -\Delta v=f, & x\in\Omega,\\ v\cdot n=0,\qquad \operatorname{curl}v\times n=0, & x\in\partial\Omega . 
	\end{cases} \]
	Thanks to \cite[Theorem~10.5]{ADN1964}, we have the following lemma, which will be used in Section 5.
	\begin{lemma}\label{lm28} 
		Let \(1<p<\infty\), \(k=0,1,\cdots\), and let \(\Omega\subset\mathbb R^N\), \(N=2,3\), be a bounded smooth domain. 
		Assume that \(v\) is a smooth solution of the above elliptic system. Then
		\begin{align}\label{x24} 
			\|v\|_{W^{k+2,p}(\Omega)} \le C \left( \|\Delta v\|_{W^{k,p}(\Omega)} + \|v\|_{L^p(\Omega)} \right). 
		\end{align} 
		In particular, when \(p=2\), 
		\begin{align}\label{x26} 
			\|v\|_{H^{k+2}(\Omega)} \le C \left( \|\Delta v\|_{H^k(\Omega)} + \|v\|_{L^2(\Omega)} \right). 
		\end{align} 
	\end{lemma}
	\pf 
	In the sense of ADN (see \cite{ADN1964}), we take 
	\( s_i=0, t_j=2, i,j=1,\cdots,N. \) After flattening the boundary, the principal boundary operators are \( v_N=0, \partial_N v_\alpha-\partial_\alpha v_N=0,  \alpha=1,\cdots,N-1. \) A direct calculation shows that these boundary operators satisfy the complementing boundary condition. Hence, by the ADN estimate \cite[Theorem~10.5]{ADN1964}, we obtain
	\[ \|v\|_{W^{k+2,p}(\Omega)} \le C \left( \|f\|_{W^{k,p}(\Omega)} + \|v\|_{L^p(\Omega)} \right). \] 
	This proves \eqref{x24}. Therefore, taking \(f=-\Delta v\) gives \eqref{x24}, and \eqref{x26} follows by choosing \(p=2\). 
	\thatsall
	
	\section{Upper bound for the density
    }
	

 This section aims to  establish the upper bound of the density. 
First, we will give out the local existence of strong solutions to the problem \eqref{x11}--\eqref{x15}, whose proof is standard in view of \cite{Kotschote,HattoriLi1994}.
\begin{lemma}\label{L41}
Under the assumptions in Theorem \ref{th11},  there exists a time $T>0$  such that the problem \eqref{x11}--\eqref{x15} admits a unique strong solution $(\rho,u)$  on $\Omega\times (0,T]$ with \eqref{x115} holds. 
\end{lemma}
    

Let $(\rho,u)$ be the local strong solution established in Lemma \ref{L41} on $Q_T=\Omega\times[0,T]$ and    $T^*>0$ be the maximal existence time, we will derive some necessary a priori estimates for any  $0<T<T^*$,  together with the De Giorgi iteration,  
and thus obtain  the  upper bound for the density. 
Different from \cite{HMZ,HuangMeng2026} considering the periodic case, we need some new ingredients here to handle  the boundary terms. 
In particular,  by virtue of the density–effective velocity bootstrap, we succeed in enlarging 
the admissible range of $\gamma$ compared with those in \cite{HMZ,HuangMeng2026}. Here, we mainly present the proof for  the 3D case with $1<\gamma<3$, where the isothermal case $\gamma=1$ can be treated by replacing the pressure potential accordingly.

Defined the initial energy as
	\[
	E_0:=
	\int_\Omega
	\left(
	\frac{1}{2}\rho_0|u_0|^2
	+
	\frac{1}{\gamma-1}\rho_0^\gamma
	+
	2|\nabla\sqrt{\rho_0}|^2
	\right)\,dx,
	\]
we will show the following the basic energy estimate.
	
	\begin{proposition}\label{pr31}
		There exists a constant \(C>0\), depending only on \(\gamma\) and \(E_0\), such that
		\begin{align}\label{x31}
			\sup_{0\le t\le T}
			\int_\Omega
			\left(
			\rho|u|^2+\rho^\gamma+|\nabla\sqrt{\rho}|^2
			\right)\,dx
			+
			\int_{Q_T}\rho|D(u)|^2\,dxdt
			\le C.
		\end{align}
	\end{proposition}
	\pf
	By \eqref{x11}\(_1\), the momentum equation \eqref{x11}\(_2\) can be rewritten as
	\[
	\rho u_t+\rho u\cdot\nabla u+\nabla\rho^\gamma
	-\operatorname{div}(2\rho D(u))
	=
	\operatorname{div}\mathbb K .
	\]
	Multiplying it by \(u\) and integrating the resulting equality over \(\Omega\), we obtain after  integration by parts that
	\begin{equation}\label{x32}
		\begin{aligned}
			&\frac{1}{2}\frac{d}{dt}\int_\Omega \rho|u|^2\,dx
			+
			\frac{1}{\gamma-1}\frac{d}{dt}\int_\Omega \rho^\gamma\,dx
			+
			2\int_\Omega \rho|D(u)|^2\,dx \\
			&~+
			\frac{1}{2}\int_{\partial\Omega}\rho|u|^2u\cdot n\,dS
			+
			\frac{\gamma}{\gamma-1}\int_{\partial\Omega}\rho^\gamma u\cdot n\,dS
			-
			2\int_{\partial\Omega}\rho(D(u)n)\cdot u\,dS\\
			&=
			\int_\Omega \operatorname{div}\mathbb K\cdot u\,dx . 
		\end{aligned}
	\end{equation}
    
	Let
	\(
	A:=\Delta\log\rho+\frac{1}{2}|\nabla\log\rho|^2,
	\)
it holds
\begin{equation}\label{x33}
		\begin{aligned}
			\int_\Omega \operatorname{div}\mathbb K\cdot u\,dx
			&=
			\int_\Omega \rho u\cdot\nabla A\,dx \\
			&=
			-\int_\Omega A\operatorname{div}(\rho u)\,dx
			+
			\int_{\partial\Omega}A\rho u\cdot n\,dS\\
			&=
			\int_\Omega A\rho_t\,dx
			+
			\int_{\partial\Omega}A\rho u\cdot n\,dS . 
		\end{aligned}
	\end{equation}
Direct calculations show that
	\begin{equation}\label{x34}
		\begin{aligned}
			\frac{d}{dt}
			\left(
			2\int_\Omega|\nabla\sqrt{\rho}|^2\,dx
			\right)
			&=
			\frac{d}{dt}
			\left(
			\frac{1}{2}\int_\Omega\frac{|\nabla\rho|^2}{\rho}\,dx
			\right) \\
			&=
			-\int_\Omega
			\left(
			\Delta\log\rho+\frac{1}{2}|\nabla\log\rho|^2
			\right)\rho_t\,dx
			+
			\int_{\partial\Omega}\frac{\rho_t}{\rho}\partial_n\rho\,dS \\
			&=
			-\int_\Omega A\rho_t\,dx
			+
			\int_{\partial\Omega}\frac{\rho_t}{\rho}\partial_n\rho\,dS .
		\end{aligned}
	\end{equation}
	Combining \eqref{x33} with \eqref{x34} leads to 
	\begin{align}\label{x35}
		\int_\Omega \operatorname{div}\mathbb K\cdot u\,dx
		=
		-\frac{d}{dt}
		\left(
		2\int_\Omega|\nabla\sqrt{\rho}|^2\,dx
		\right)
		+
		\int_{\partial\Omega}\frac{\rho_t}{\rho}\partial_n\rho\,dS
		+
		\int_{\partial\Omega}A\rho u\cdot n\,dS .
	\end{align}
	Substituting \eqref{x35} into \eqref{x32} yields
	\begin{align}\label{x36}
		&\frac{d}{dt}
		\int_\Omega
		\left(
		\frac{1}{2}\rho|u|^2
		+
		\frac{1}{\gamma-1}\rho^\gamma
		+
		2|\nabla\sqrt{\rho}|^2
		\right)\,dx
		+
		2\int_\Omega \rho|D(u)|^2\,dx
		+
		\mathcal B_{u}(t)
		=0,
	\end{align}
	where 
	\begin{equation}\label{x37}
		\begin{aligned}
			\mathcal B_{u}(t)
			&=
			\frac{1}{2}\int_{\partial\Omega}\rho|u|^2u\cdot n\,dS
			+
			\frac{\gamma}{\gamma-1}\int_{\partial\Omega}\rho^\gamma u\cdot n\,dS
			-
			2\int_{\partial\Omega}\rho(D(u)n)\cdot u\,dS\\
			&\quad
			-
			\int_{\partial\Omega}\frac{\rho_t}{\rho}\partial_n\rho\,dS
			-
			\int_{\partial\Omega}A\rho u\cdot n\,dS\\
            &=-2\int_{\partial\Omega}\rho(D(u)n)\cdot u\,dS
		\end{aligned}
	\end{equation}
	owing to the boundary condition \eqref{x15}.
    Moreover, one has
	\[
	(D(u)n)\cdot u
	=
	(D(u)n)_{tan}\cdot u_{tan}
	=
	-\kappa(x)u_{tan}\cdot u_{tan}
	\quad\text{on }\partial\Omega,
	\]
	which along with \eqref{x19} gives
	\begin{align}\label{x38}
		\mathcal B_{u}(t)
		=
		2\int_{\partial\Omega}\rho\,\kappa(x)u_{tan}\cdot u_{tan}\,dS
		\ge0.
	\end{align}
	This combined with  \eqref{x36} deduces that
	\[
	\frac{d}{dt}
	\int_\Omega
	\left(
	\frac{1}{2}\rho|u|^2
	+
	\frac{1}{\gamma-1}\rho^\gamma
	+
	2|\nabla\sqrt{\rho}|^2
	\right)\,dx
	+
	2\int_\Omega \rho|D(u)|^2\,dx
	\le0.
	\]
	Integrating it over \([0,T]\) yields
	\[
	\sup_{0\le t\le T}
	\int_\Omega
	\left(
	\rho|u|^2+\rho^\gamma+|\nabla\sqrt{\rho}|^2
	\right)\,dx
	+
	\int_{Q_T}\rho|D(u)|^2\,dxdt
	\le C(\gamma,E_0),
	\]
	and thus completes the proof of Proposition~\ref{pr31}.
	\thatsall\\
    
    Next, we derive the corresponding energy estimate for the effective velocity $v$.
	\begin{proposition}\label{pr32}
		There exists a constant \(C>0\), depending only on \(\gamma\) and \(E_0\), such that
		\begin{align}\label{x39}
			\sup_{0\le t\le T}\int_\Omega \rho|v|^2\,dx
			+
			\int_{Q_T}|\nabla\rho^{\frac{\gamma}{2}}|^2\,dxdt
			+
			\int_{Q_T}\rho|\nabla v|^2\,dxdt
			\le C.
		\end{align}
	\end{proposition}
	\pf
    Integrating \eqref{x16}\(_2\) multiplied by \(\rho v\)  over \(\Omega\),  one derives after  using \eqref{x16}\(_1\) and integration by parts   that
	\begin{equation}\label{x310}
		\begin{aligned}
			&\int_\Omega \rho v_t\cdot v\,dx
			+
			\int_\Omega \rho(v-\nabla\log\rho)\cdot\nabla v\cdot v\,dx\\
			&+
			\int_\Omega \nabla\rho^\gamma\cdot v\,dx
			-
			\int_\Omega \operatorname{div}(\rho\nabla v)\cdot v\,dx \\
			=&
			\frac{1}{2}\frac{d}{dt}\int_\Omega \rho|v|^2\,dx
			+
			\frac{1}{\gamma-1}\frac{d}{dt}\int_\Omega \rho^\gamma\,dx
			+
			\frac{4}{\gamma}\int_\Omega |\nabla\rho^{\frac{\gamma}{2}}|^2\,dx\\
			&+
			\int_\Omega \rho|\nabla v|^2\,dx
			+ 
			\mathcal B_v(t)\\
			=&
			0,
		\end{aligned}
	\end{equation}
	in which
    \begin{equation}\label{x311}
		\begin{aligned}
			\mathcal B_v(t)
			&=
			\frac{1}{2}\int_{\partial\Omega}\rho|v|^2v\cdot n\,dS
			-
			\frac{1}{2}\int_{\partial\Omega}|v|^2\partial_n\rho\,dS\\
			&\quad+
			\frac{\gamma}{\gamma-1}\int_{\partial\Omega}\rho^\gamma u\cdot n\,dS
			-
			\int_{\partial\Omega}\rho(n\cdot\nabla v)\cdot v\,dS\\
            &=-\int_{\partial\Omega}\rho(n\cdot\nabla v)\cdot v\,dS
		\end{aligned}
	\end{equation}
	due to the boundary conditions \eqref{x15} and \eqref{x18}.
   
    We first identify the normal derivative of \(v\) on the boundary.
	For any tangential vector \(\eta\in T_x(\partial\Omega)\), it follows from the boundary
	conditions \eqref{x18} that 
	\[
	(n\cdot\nabla v)\cdot\eta
	=
	(\eta\cdot\nabla v)\cdot n.
	\]
	Furthermore, applying $\nabla$ to  \(v\cdot n=0\)   on the direction \(\eta\)  gives
	\[
	(\eta\cdot\nabla v)\cdot n
	=
	-(\eta\cdot\nabla n)\cdot v.
	\]
	Notice that
	\(
	\kappa(x)\eta=(\eta\cdot\nabla)n,
	\)
	and recalling that \(\kappa(x)\) is symmetric on \(T_x(\partial\Omega)\), we obtain
	\begin{align}\label{x3111}
		(n\cdot\nabla v)\cdot\eta
		=
		(\eta\cdot\nabla v)\cdot n 
		=
		-(\eta\cdot\nabla n)\cdot v 
		=
		-\kappa(x)v_{tan}\cdot\eta .
	\end{align}
	Since $$v=v_{tan},~~~\text{on} ~~\partial\Omega,$$ choosing
	\(\eta=v_{tan}\) in \eqref{x3111} gives
	\begin{align}\label{x3121}
		(n\cdot\nabla v)\cdot v
		=
		-\kappa(x)v_{tan}\cdot v_{tan} .
	\end{align}
	Consequently, we have
	\begin{align}\label{x312}
		\mathcal B_v(t)
		=
		\int_{\partial\Omega}
		\rho\,\kappa(x)v_{tan}\cdot v_{tan}\,dS
		\ge 0.
	\end{align}
	Hence, one deduces from \eqref{x310} and \eqref{x312}  that
	\begin{equation}\label{x313}
		\begin{aligned}
			&\sup_{0\le t\le T}
			\int_\Omega
			\left(
			\rho|v|^2+\rho^\gamma
			\right)\,dx
			+
			\int_{Q_T}|\nabla\rho^{\frac{\gamma}{2}}|^2\,dxdt
			+
			\int_{Q_T}\rho|\nabla v|^2\,dxdt \\
			&\le
			C
			\int_\Omega
			\left(
			\rho_0|v_0|^2+\rho_0^\gamma
			\right)\,dx\\
            &\le C(\gamma,E_0),
		\end{aligned}
	\end{equation}
    where the last inequality follows from the definition of the initial effective velocity $v_0$ in \eqref{x17} and the initial energy $E_0$.
This completes the proof of Proposition \ref{pr32}.
	\thatsall\\

    To prove Proposition~\ref{pr33} for $1<\gamma<3$, we first establish some necessary estimates for the bootstrap argument, see Corollaries~\ref{cor31}--\ref{cor32}. We begin with the following estimates for $\sqrt{\rho}$ and $\nabla^2\log\rho$.
	\begin{corollary}\label{cor31}
		Under the assumptions of Propositions~\ref{pr31}--\ref{pr32}, there exists a constant \(C_{1}>0\), depending only on \(\gamma\) and \(E_0\), such that
		\begin{align}\label{x313a}
			\int_{Q_T}\rho|\nabla^2\log\rho|^2\,dxdt\le C_{1},
		\end{align}
		and
		\begin{align}\label{x313b}
			\int_0^T\|\nabla^2\sqrt\rho\|_{L^2}^2\,dt\le C_{1}.
		\end{align}
		Moreover, there exists a constant \(C_{2}>0\), depending on \(T^*,\gamma,\Omega\) and \(E_0\), such that
		\begin{align}\label{x313c}
			\sup_{0\le t\le T}\|\sqrt\rho(t)\|_{H^1}^2
			+
			\int_0^T
			\left(
			\|\sqrt\rho\|_{H^2}^2
			+
			\|\nabla\sqrt\rho\|_{L^6}^2
			\right)\,dt
			\le C_{2}.
		\end{align}
	\end{corollary}
	\pf
	Since
	\(
	v=u+\nabla\log\rho,
	\)
	we have
	\(
	\nabla^2\log\rho=D(v)-D(u).
	\)
	Combining this identity with \eqref{x31} and \eqref{x39}, we obtain \eqref{x313a} for some $C_1=C(\gamma,E_0)>0$. It remains to establish \eqref{x313b} and \eqref{x313c}.
	
	To prove \eqref{x313b}, we set \(f=\sqrt{\rho}\). The boundary condition \(\partial_n\rho=0\) implies that
	\[
	\partial_n f=0 \qquad\text{on }\partial\Omega.
	\]
	Using the identity
	\(
	f\nabla^2\log f
	=
	\nabla^2f-\frac{\nabla f\otimes\nabla f}{f},
	\)
	we have
	\[
	\begin{aligned}
		\int_\Omega f^2|\nabla^2\log f|^2\,dx
		&=
		\int_\Omega|\nabla^2f|^2\,dx
		+\int_\Omega\frac{|\nabla f|^4}{f^2}\,dx  
		-2\int_\Omega
		\nabla^2f:\frac{\nabla f\otimes\nabla f}{f}\,dx\\
        &=
		\int_\Omega|\nabla^2f|^2\,dx
		+\int_\Omega\frac{|\nabla f|^4}{f^2}\,dx  
		+\int_\Omega
		|\nabla f|^2\Delta\log f\,dx.
	\end{aligned}
	\]
	By H\"older's and Young's inequalities, together with
	\(
	|\Delta\log f|^2\le 3|\nabla^2\log f|^2,
	\)
	we obtain
	\[
	\int_\Omega
		|\nabla f|^2\Delta\log f\,dx
	\le
	\frac12\int_\Omega\frac{|\nabla f|^4}{f^2}\,dx
	+\frac32\int_\Omega f^2|\nabla^2\log f|^2\,dx.
	\]
	Substituting this estimate into the preceding identity yields
	\[
	\int_\Omega|\nabla^2f|^2\,dx
	\le
	\frac52\int_\Omega f^2|\nabla^2\log f|^2\,dx.
	\]
	Since \(\log\rho=2\log f\) and \(f^2=\rho\), \eqref{x313a} yields \eqref{x313b}, after enlarging $C_1$ if necessary. 
    
    Finally, by \eqref{x11}\(_1\), Proposition~\ref{pr31} and the Sobolev embedding
	\(H^2(\Omega)\hookrightarrow W^{1,6}(\Omega)\), we obtain \eqref{x313c} with $C_2=C(T^*,\gamma,\Omega,E_0)>0$.
	\thatsall\\

    Next, we need the following $L^2(0,T;L^\infty)$-estimate for $\rho$, which will be used to derive the weighted estimates for $v$ and improve the integrability of the density.
	\begin{corollary}\label{cor32}
		Assume that \(N=3\). Under the assumptions of Corollary~\ref{cor31}, there exists a constant \(C >0\), depending on \(T^*,\gamma,\Omega\) and \(E_0\), such that
		\begin{align}\label{x313d}
			\int_0^T\|\rho(t)\|_{L^\infty}^2\,dt
			\le C,
		\end{align}
	\end{corollary}
	\pf
	Set \(f=\sqrt\rho\). By \eqref{x313c} and the Sobolev embedding
	\(H^1(\Omega)\hookrightarrow L^6(\Omega)\), we have
	\[
	\sup_{0\le t\le T}\|f(t)\|_{L^6}\le C,
	\qquad
	\int_0^T\|\nabla f\|_{L^6}^2\,dt\le C.
	\]
	Applying Lemma~\ref{lm24} with \(q=r=6\), we obtain
	\[
	\|f\|_{L^\infty}
	\le
	C\|f\|_{L^6}^{1/2}\|\nabla f\|_{L^6}^{1/2}
	+C\|f\|_{L^6}.
	\]
	Since \(\rho=f^2\), it follows that
	\[
	\|\rho\|_{L^\infty}^2
	=\|f\|_{L^\infty}^4
	\le
	C\|f\|_{L^6}^2\|\nabla f\|_{L^6}^2
	+C\|f\|_{L^6}^4.
	\]
	Integrating over \([0,T]\) and using \eqref{x313c}, together with
	\(T<T^*\), gives
	\[
	\int_0^T\|\rho(t)\|_{L^\infty}^2\,dt
	\le C(T^*,\gamma,\Omega,E_0),
	\]
	which proves \eqref{x313d}.
	\thatsall\\

    Now, we use a bootstrap argument to establish the following weighted $L^{q+2}$-estimate for  $v$, which plays a crucial role in the subsequent De Giorgi iteration for the upper bound of the density.
	\begin{proposition}\label{pr33}
		Assume that \(N=3\) and \(1<\gamma<3\). Let \(q>1\) be fixed. Then there exists a constant \(C>0\), depending on \(T^*,q,\gamma,E_0,\Omega,\ \|\rho_0\|_{L^\infty} \)and  \(\|\rho_0^{1/(q+2)}v_0\|_{L^{q+2}},\) such that
		\begin{align}\label{x315}
			\sup_{0\le t\le T}
			\left\|\rho^{\frac{1}{q+2}}v\right\|_{L^{q+2}}
			\le C.
		\end{align}
	\end{proposition}
	\pf
	The proof is divided into two steps. We first derive two estimates which improve the integrability of $v$ and $\rho$, respectively, and then use a bootstrap argument to complete the proof.
	
	\medskip
	\noindent
	{\bf Step 1.}
	Let \(Q\ge3\) and suppose that
	\[
	\sup_{0\le t\le T}\int_\Omega\rho^Q\,dx\le C.
	\]
	For \(p>2\), multiplying \eqref{x16}\(_2\) by
	\(\rho|v|^{p-2}v\), integrating over \(\Omega\), and using
	\eqref{x16}\(_1\), we obtain
	\[
	\begin{aligned}
		&\frac{1}{p}\frac{d}{dt}\int_\Omega\rho|v|^p\,dx
		+\int_\Omega\rho|v|^{p-2}|\nabla v|^2\,dx
		+(p-2)\int_\Omega\rho|v|^{p-2}|\nabla|v||^2\,dx
		+\mathcal B_{\rho}(t)\\
		&\qquad
		=\int_\Omega\rho^\gamma
		\operatorname{div}(|v|^{p-2}v)\,dx .
	\end{aligned}
	\]
	By \eqref{x18} and \eqref{x3121}, the only boundary term is
	\[
	\mathcal B_{\rho}(t)
	=
	\int_{\partial\Omega}
	\rho|v|^{p-2}\kappa(x)v_{tan}\cdot v_{tan}\,dS
	\ge0,
	\]
	while all the other boundary terms vanish due to
	\(v\cdot n=0\) and \(\partial_n\rho=0\).
	Young's and H\"older's inequalities then give
	\[
	\frac{d}{dt}
	\left\|\rho^{1/p}v\right\|_{L^p}^2
	\le
	C
	\left(
	\int_\Omega\rho^{(\gamma-1)p+1}\,dx
	\right)^{\frac{2}{p}}.
	\]
	If \((\gamma-1)p+1\le Q\), the right-hand side is bounded by
	the assumed \(L^Q\)-estimate of \(\rho\). Otherwise,
	\[
	\left(
	\int_\Omega\rho^{(\gamma-1)p+1}\,dx
	\right)^{\frac{2}{p}}
	\le
	C
	\|\rho\|_{L^\infty}^{
		\frac{2((\gamma-1)p+1-Q)}{p}}.
	\]
    Hence, under the condition
    \begin{align}\label{x320a}
	(\gamma-2)p<Q-1,
    \end{align}
    Corollary~\ref{cor32} and H\"older's inequality in time yields
    \begin{align}\label{x321}
	\sup_{0\le t\le T}
	\left\|\rho^{1/p}v\right\|_{L^p}
	\le
	C\left(
	1+\left\|\rho_0^{1/p}v_0\right\|_{L^p}
	\right).
    \end{align}
	
	Next, using the weighted estimate \eqref{x321}, we improve the integrability of the density. Let $R>Q$, multiplying \eqref{x16}$_1$ by $\rho^{R-1}$ and integrating over $\Omega$ by parts, we obtain
	\[
	\frac{1}{R}\frac{d}{dt}\int_\Omega\rho^R\,dx
	+(R-1)\int_\Omega\rho^{R-2}|\nabla\rho|^2\,dx
	=
	(R-1)\int_\Omega\rho^{R-1}v\cdot\nabla\rho\,dx.
	\]
	Here, the boundary contributions vanish exactly since \eqref{x18}.
	Hence, by Young's inequality, a direct calculation shows that
	\[
	\frac{d}{dt}\int_\Omega\rho^R\,dx
	+\frac{1}{2}R(R-1)\int_\Omega\rho^{R-2}|\nabla\rho|^2\,dx
	\le
	C\int_\Omega\rho^R|v|^2\,dx.
	\]
	Using \eqref{x321} and H\"older's inequality, we have
	\[
	\begin{aligned}
		\int_\Omega\rho^R|v|^2\,dx
		&\le
		\left\|\rho^{1/p}v\right\|_{L^p}^2
		\left(
		\int_\Omega
		\rho^{\frac{pR-2}{p-2}}\,dx
		\right)^{\frac{p-2}{p}}\\
		&\le
		C
		\|\rho\|_{L^\infty}^{
			R-Q+\frac{2(Q-1)}{p}},
	\end{aligned}
	\]
	where the last inequality follows from \(R>Q\) and the \(L^Q\)-bound of \(\rho\). 
    
    Thus, combining the weighted estimate \eqref{x321} with the assumed $L^Q$-bound for $\rho$, we see that if
    \begin{align}\label{x321b}
    R-Q+\frac{2(Q-1)}{p}\le 2,
    \end{align} then Corollary~\ref{cor32} yields the improved $L^R$-estimate for the density:\begin{align}\label{x321c}
    \sup_{0\le t\le T}\int_\Omega\rho^R dx\le C.\end{align}
	
	\medskip
	\noindent
	{\bf Step 2.}
	Employing Corollary~\ref{cor31} yields
	\[
	\sup_{0\le t\le T}\|\rho(t)\|_{L^3}
	=
	\sup_{0\le t\le T}\|\sqrt\rho(t)\|_{L^6}^2
	\le C.
	\]
	Let
	\(
	p_*:=q+2>3.
	\)
    
	For $1<\gamma\le2$, we have $(\gamma-2)p_*<2=3-1$. Hence, taking $Q=3$ and $p=p_*$ in \eqref{x321}, we directly obtain \eqref{x315}.
    
    It remains to consider the case \(2<\gamma<3\), define
	\[
	Q_0:=3,\qquad
	Q_{k+1}:=Q_k+3-\gamma,
	\qquad
	p_k:=\frac{2(Q_k-1)}{\gamma-1}.
	\]
	Choose \(N\) to be the smallest nonnegative integer such that
	\[
	(\gamma-2)p_*<Q_N-1.
	\]
    For every \(0\le k<N\), the minimality of \(N\) implies
	\(
	Q_k-1\le(\gamma-2)p_*.
	\)
	It follows that
	\[
	2<p_k
	\le
	\frac{2(\gamma-2)}{\gamma-1}p_*
	<p_*,\qquad
	(\gamma-2)p_k
	=
	\frac{2(\gamma-2)}{\gamma-1}(Q_k-1)
	<Q_k-1.
	\]
	Moreover, by the choice of $p_k$ and $Q_{k+1}$, we have
	\[
	Q_k+2-\frac{2(Q_k-1)}{p_k}
	=
	Q_k+3-\gamma
	=
	Q_{k+1}.
	\]
	Therefore, conditions \eqref{x320a} and \eqref{x321b} are satisfied with $Q=Q_k$, $p=p_k$ and $R=Q_{k+1}$. Applying \eqref{x321} and \eqref{x321c} successively, and starting from $Q_0=3$, we obtain
	\[
	\sup_{0\le t\le T}\int_\Omega\rho^{Q_k}\,dx\le C,
	\qquad
	k=0,\ldots,N.
	\]
    By the definition of \(N\), applying \eqref{x320a} once more with \(Q=Q_N\) and \(p=p_*\) gives
	\[
	\sup_{0\le t\le T}
	\left\|\rho^{1/p_*}v\right\|_{L^{p_*}}
	\le C.
	\]
	Since \(p_*=q+2\), this proves \eqref{x315}.
	\thatsall\\
	
	We now turn to the De Giorgi iteration for the density equation. For \(k\ge2\), we denote that
	\[
	\rho^{(k)}:=(\rho-k)^+,\qquad
	\Gamma_k(t):=\{x\in\Omega:\rho(x,t)\ge k\},
	\]
	together with
	\[
	|\rho^{(k)}|_{Q_t}
	:=
	\|\rho^{(k)}\|_{L_t^\infty L^2}
	+
	\|\nabla\rho^{(k)}\|_{L_t^2L^2},\qquad
	\mu(k):=\int_0^T |\Gamma_k(t)|^{\frac{q}{q+2}}\,dt .
	\]
	Here and below, \(|\Gamma_k(t)|\) denotes the Lebesgue measure of \(\Gamma_k(t)\).

    We begin with the following estimate for the truncated density, using the method of \cite[Lemma~2.1]{HMZ} and taking the boundary terms into account.
	\begin{lemma}\label{lm31}
		There exists a constant \(C>0\), depending on
		\(T^*,q,\gamma,E_0,\Omega,\ \|\rho_0\|_{L^\infty} \)and  \(\|\rho_0^{1/(q+2)}v_0\|_{L^{q+2}},\)such that for any \(k\ge2\),
		\begin{align}\label{x322}
			|\rho^{(k)}|_{Q_T}
			\le
			Ck\mu(k)^{\frac{1}{2}}
			+
			C\|\rho_0^{(k)}\|_{L^2}.
		\end{align}
		In particular, if \(k\ge \|\rho_0\|_{L^\infty}\), then
		\[
		|\rho^{(k)}|_{Q_T}
		\le
		Ck\mu(k)^{\frac{1}{2}}.
		\]
	\end{lemma}
	\pf
	Since \(k\) is a constant, the first equation of \eqref{x16}\(_1\) gives
	\[
	\rho_t^{(k)}+\operatorname{div}(\rho v)-\Delta\rho^{(k)}=0
	\quad\text{on }\Gamma_k(t).
	\]
	Testing it by \(\rho^{(k)}\), integrating over \(\Omega\), and using \eqref{x18}, we get
	\begin{equation}\label{x323}
		\begin{aligned}
			\frac{1}{2}\frac{d}{dt}\int_\Omega |\rho^{(k)}|^2\,dx
			+
			\int_\Omega |\nabla\rho^{(k)}|^2\,dx
			&=
			\int_{\Gamma_k(t)}\rho v\cdot\nabla\rho^{(k)}\,dx \\
			&\le
			\frac{1}{2}\int_\Omega |\nabla\rho^{(k)}|^2\,dx
			+
			C\int_{\Gamma_k(t)}\rho^2|v|^2\,dx .
		\end{aligned}
	\end{equation}
	Hence, for some constant \(C_3>0\), a direct calculation shows that
	\begin{align}\label{x324}
		|\rho^{(k)}|_{Q_T}^2
		\le
		C_3\int_0^T\int_{\Gamma_k(t)}\rho^2|v|^2\,dxdt
		+
		C_3\|\rho_0^{(k)}\|_{L^2}^2 .
	\end{align}
	
	It remains to estimate the first term on the right-hand side. By
	Lemma~\ref{lm21}, Proposition~\ref{pr33}, Young's and
	H\"older's inequality, for any \(\varepsilon>0\), we have
	\begin{equation}\label{x325}
		\int_0^T\int_{\Gamma_k(t)}\rho^2|v|^2\,dxdt
		\le
		\varepsilon|\rho^{(k)}|_{Q_T}^2
		+
		C(\varepsilon) k^2\mu(k).
	\end{equation}
	Choosing \(\varepsilon>0\) sufficiently small such that
	\(C_3\varepsilon\le\frac12\), and combining \eqref{x324} with
	\eqref{x325}, we obtain
	\[
	|\rho^{(k)}|_{Q_T}^2
	\le
	Ck^2\mu(k)
	+
	C\|\rho_0^{(k)}\|_{L^2}^2,
	\]
	which implies \eqref{x322}.
	\thatsall\\
	
	To apply the preceding estimate in the De Giorgi iteration, we need the following lemma from \cite[Lemma 2.2]{YuWu}.
	
	\begin{lemma}\label{lm32}
		Let \(\{\Phi_n\}_{n\ge0}\) be a nonnegative sequence satisfying
		\[
		\Phi_{n+1}\le c b^n\Phi_n^{1+\eta},
		\qquad n=0,1,2,\cdots,
		\]
		where \(c>0\), \(b>1\), and \(\eta>0\). If
		\[
		\Phi_0\le c^{-1/\eta}b^{-1/\eta^2},
		\]
		then
		\[
		\Phi_n\le c^{-1/\eta}b^{-1/\eta^2}b^{-n/\eta},
		\qquad
		\lim_{n\to\infty}\Phi_n=0.
		\]
	\end{lemma}
	
	Finally, we establish the upper bound for the density by combining Lemma~\ref{lm31}--\ref{lm32}.
	
	\begin{proposition}\label{pr34}
		There exists a constant \(C>0\), depending on
		\(T^*,q,\gamma,E_0,\Omega,\|\rho_0^{1/(q+2)}v_0\|_{L^{q+2}}\), and
		\(\|\rho_0\|_{L^\infty}\), such that
		\begin{align}\label{x328}
			\sup_{0\le t\le T}\|\rho(t)\|_{L^\infty}\le C.
		\end{align}
	\end{proposition}
	\pf
	We denote that
	\[
	\widehat K_0:=\max\{2,\|\rho_0\|_{L^\infty}\},
	\qquad
	K_n:=M\widehat K_0(2-2^{-n}),
	\qquad n=0,1,2,\cdots,
	\]
	where \(M>1\) will be chosen later. A direct calculation shows that
	\begin{equation}\label{x329}
		\begin{aligned}
			\|\rho^{(K_n)}\|_{L_T^{\frac{10q+8}{3q+6}}L^{\frac{10q+8}{3q}}}
			\le C|\rho^{(K_n)}|_{Q_T}.
		\end{aligned}
	\end{equation}
	Using Lemma~\ref{lm23}, \eqref{x329}, together with Lemma~\ref{lm31}, we have
	\begin{equation}\label{x330}
		\begin{aligned}
			\mu(K_{n+1})
			\le
			C
			\left(
			\frac{K_n}{K_{n+1}-K_n}
			\right)^{\frac{10q+8}{3q+6}}
			\mu(K_n)^{\frac{5q+4}{3q+6}} .
		\end{aligned}
	\end{equation}
	Thus, by the definition of $K_n$ and denoting the constant $C$ by $C_4$, we obtain the relation
    \begin{align}\label{x331}
    \mu(K_{n+1})
    \le
    C_4\,16^n\,\mu(K_n)^{1+\eta}, \qquad
    \eta:=\frac{2q-2}{3q+6}>0.
    \end{align}
	
	It remains to verify the smallness condition required by Lemma~\ref{lm32}. Using \eqref{x11}\(_1\) and \(K_0=M\widehat K_0\), we have
	\begin{equation}\label{x332}
		\begin{aligned}
			\mu(K_0)
			\le
			\int_0^T
			\left(
			\frac{1}{M\widehat K_0}\int_\Omega \rho(x,t)\,dx
			\right)^{\frac{q}{q+2}}dt
			\le C_5 M^{-\frac{q}{q+2}},
		\end{aligned}
	\end{equation}
    where $C_5>0$ is independent of $M$.
	We now choose \(M>1\) sufficiently large so that
	\begin{align}\label{x333}
		C_5M^{-\frac{q}{q+2}}
		\le
		C_4^{-1/\eta}16^{-1/\eta^2}.
	\end{align}
	Then, applying Lemma~\ref{lm32} to \eqref{x331}, we obtain
	which implies that
	\[
	\lim_{n\to\infty}\mu(K_n)=0.
	\]
	By the continuity from above of the Lebesgue measure, we have $\mu(2M\widehat K_0)=0$, which implies
	\[
    |\{x\in\Omega:\rho(x,t)\ge 2M\widehat K_0\}|=0,\qquad \text{for a.e. }t\in[0,T].
    \] 
    Since $\rho\in C([0,T];H^3(\Omega))\hookrightarrow C([0,T];L^\infty(\Omega))$, the above estimate holds for all $t\in[0,T]$. Consequently,
    \[
    \sup_{0\le t\le T}\|\rho(t)\|_{L^\infty}
    \le 2M\widehat K_0.
    \]
    This completes the proof of Proposition~\ref{pr34}.
	\thatsall
	
	\section{Positive lower bound for the density (\(N=2,3\))}
	
	In this section, we establish the positive lower bound of the density, which is a crucial step in extending the local strong solution globally in time.
	
	In contrast to the periodic case, the boundary terms arising from integration by parts require additional treatment. Here, we use boundary condition \eqref{x18}, and the assumption \eqref{x19} to control these terms. We begin with Lemma~4.1, which establishes a weighted $L^{p+2}$-estimate for $v$.
	
	\begin{lemma}\label{lm41}
		Let \(p>2\). There exists a constant \(\widehat C\ge1\), independent of \(p\), such that
		\begin{align}\label{x41}
			\sup_{0\le t\le T}
			\|\rho^{\frac{1}{p+2}}v\|_{L^{p+2}}
			\le
			\widehat C\sqrt{p+2}.
		\end{align}
		Here \(\widehat C\) may depend on
		\(T^*,q,\gamma,E_0,\Omega,\|\rho_0\|_{L^\infty},
		\|v_0\|_{L^\infty}\), and
		\(\|\rho_0^{1/(q+2)}v_0\|_{L^{q+2}}\), but not on \(p\).
	\end{lemma}
	\pf
	Multiplying the second equation of \eqref{x16}\(_2\) by \(\rho|v|^pv\), and using \eqref{x16}\(_1\), we obtain
	\begin{equation}\label{x42}
		\begin{aligned}
			&\frac{1}{p+2}\frac{d}{dt}\int_\Omega \rho|v|^{p+2}\,dx
			+
			\int_\Omega \rho|v|^p|\nabla v|^2\,dx
			+
			p\int_\Omega \rho|v|^p|\nabla |v||^2\,dx
			+
			\mathcal B_{\rho}(t)\\
			=&
			\int_\Omega \rho^\gamma\operatorname{div}(|v|^pv)\,dx ,
		\end{aligned}
	\end{equation}
	where the boundary contribution is given by
	\begin{equation}\label{x43}
		\begin{aligned}
			\mathcal B_{\rho}(t)
			&=
			\frac{1}{p+2}\int_{\partial\Omega}\rho|v|^{p+2}v\cdot n
			-\frac{1}{p+2}\int_{\partial\Omega}|v|^{p+2}\partial_n\rho
			+\int_{\partial\Omega}\rho^\gamma |v|^pv\cdot n\\
			&\quad
			-\int_{\partial\Omega}\rho |v|^p(n\cdot\nabla v)\cdot v .
		\end{aligned}
	\end{equation}
	Using \eqref{x18}, the first three terms in \eqref{x43} vanish. For the remaining term, \eqref{x3121} gives
    \[
    (n\cdot\nabla v)\cdot v
    =
    -\kappa(x)v_{tan}\cdot v_{tan}
    \quad\text{on }\partial\Omega.
    \]
	Hence
	\begin{align}\label{x44}
		\mathcal B_\rho(t)
		=
		\int_{\partial\Omega}\rho |v|^p\kappa(x)v_{tan}\cdot v_{tan}
		\ge0.
	\end{align}
	
	Returning to \eqref{x42}, we may discard the nonnegative boundary term. Proposition~\ref{pr32}, Proposition~\ref{pr34}, together with H\"older's and Young's inequalities yields
	\[
	\frac{d}{dt}{\int_\Omega \rho|v|^{p+2}\,dx}^{\frac{2}{p+2}}
	\le
	C(p+1).
	\]
	Integrating over \([0,t]\subset[0,T]\) gives
	\begin{align}\label{x47}
		\|\rho^{\frac{1}{p+2}}v(t)\|_{L^{p+2}}^{2}
		\le
		\|\rho_0^{\frac{1}{p+2}}v_0\|_{L^{p+2}}^{2}
		+
		C(p+1).
	\end{align}
    where the dependence on \(T\) is absorbed into the constant \(C\).
    
	It remains to estimate the initial term. Since
	\(\rho_0\in L^\infty(\Omega)\), \(v_0\in L^\infty(\Omega)\), and \(\Omega\) is bounded, a direct calculation shows that
	\[
	\|\rho_0^{\frac{1}{p+2}}v_0\|_{L^{p+2}}
	\le C
	\]
	uniformly for \(p>2\). Estimate \eqref{x47} therefore implies
	\[
	\sup_{0\le t\le T}
	\|\rho^{\frac{1}{p+2}}v\|_{L^{p+2}}
	\le
	\widehat C\sqrt{p+2},
	\]
	which proves \eqref{x41}.
	\thatsall
	
	\begin{remark}\label{rm41}
		Estimate \eqref{x41} grows like \(\sqrt{p}\) as \(p\to\infty\). Hence it does not provide a uniform \(L^p\)-control of \(v\) with respect to \(p\). This is the obstruction to applying directly the lower-bound argument based on a uniform-in-\(p\) estimate of the effective velocity.
	\end{remark}
	
	The following proposition gives the key control relation between the effective velocity and the lower bound of the density. We set
	\[
	\tau:=\rho^{-1},
	\qquad
	\Theta_T:=
	\sup_{0\le t\le T}\|\tau(t)\|_{L^\infty}
	+
	e^{\frac{25}{9}}.
	\]
	
	\begin{proposition}\label{pr42}
		There exists a constant \(C_v\ge 1\), depending on
		\(T^*,q,\gamma,E_0,\Omega,\)\(\|\rho_0\|_{L^\infty},\)
		\(\|v_0\|_{L^\infty}\), and \(\|\rho_0^{\frac{1}{q+2}}v_0\|_{L^{q+2}}\), such that
		\begin{align}\label{x48}
			\sup_{0\le t\le T}\|v(t)\|_{L^\infty}
			\le
			C_v\sqrt{\log\Theta_T}.
		\end{align}
	\end{proposition}
	\pf
	Let \(p>2\). Multiplying \eqref{x16}\(_2\) by \(\rho|v|^pv\), integrating over \(\Omega\), and arguing exactly as in Proposition~\ref{pr33} with \(q\) replaced by \(p\), the boundary contribution is nonnegative and hence can be discarded. It follows that
	\begin{equation}\label{x49}
		\begin{aligned}
			&\frac{1}{p+2}\frac{d}{dt}\int_\Omega \rho|v|^{p+2}\,dx
			+
			\frac{1}{2}\int_\Omega \rho|v|^p|\nabla v|^2\,dx
			+
			\frac{p}{2}\int_\Omega \rho|v|^p|\nabla |v||^2\,dx\\
			\le&C(p+1)\int_\Omega \rho^{2\gamma-1}|v|^p\,dx .
		\end{aligned}
	\end{equation}
	After integrating \eqref{x49} over \([0,T]\), Proposition~\ref{pr34}, Proposition~\ref{pr32}, and H\"older's inequality in time gives
	\begin{equation}\label{x410}
		\begin{aligned}
			&\sup_{0\le t\le T}
			\frac{1}{p+2}\int_\Omega \rho|v|^{p+2}\,dx
			+
			\frac{1}{2}\int_{Q_T}\rho|v|^p|\nabla v|^2\,dxdt
			+
			\frac{p}{2}\int_{Q_T}\rho|v|^p|\nabla |v||^2\,dxdt\\
			\le&
			C(p+1)
			\left(\int_{Q_T}\rho|v|^{p+2}\,dxdt\right)^{\frac{p-2}{p}}
			+
			\frac{2}{p+2}\int_\Omega \rho_0|v_0|^{p+2}\,dx .
		\end{aligned}
	\end{equation}
	
	We now improve the space--time integrability of \(v\). In view of the upper bound for \(\rho\) and the fact that \(\tau\le\Theta_T\), H\"older's inequality, Lemma~\ref{lm24}, and \eqref{x410}, we obtain the following estimate, where we may assume $C_6\ge1$ without loss of generality:
    \begin{equation}\label{x411}
		\begin{aligned}
			\int_{Q_T}\rho|v|^{\frac{5}{3}(p+2)}\,dxdt
			&\le
			\int_0^T
			\left(\int_\Omega \rho^{\frac{3}{2}}|v|^{p+2}\,dx\right)^{\frac{2}{3}}
			\left(\int_\Omega |v|^{3(p+2)}\,dx\right)^{\frac{1}{3}}dt\\
			&\le
			C\Theta_T
			\left(
			\sup_{0\le t\le T}\int_\Omega \rho|v|^{p+2}\,dx
			+
			\int_{Q_T}\rho
			\left|\nabla |v|^{\frac{p+2}{2}}\right|^2\,dxdt
			\right)^{\frac{5}{3}}\\
			&\le
			C_6\Theta_T(p+2)^{\frac{10}{3}}
			\left(\int_{Q_T}\rho|v|^{p+2}\,dxdt\right)^{\frac{5}{3}}
			+
			C_6\Theta_T(p+2)^{\frac{10}{3}}\\
			&\quad
			+
			C_6\Theta_T
			\left(\int_\Omega \rho_0|v_0|^{p+2}\,dx\right)^{\frac{5}{3}}.
		\end{aligned}
	\end{equation}
    To estimate the initial term in \eqref{x411}, we apply H\"older's and Young's inequalities to obtain
	\begin{equation}\label{x412}
		\begin{aligned}
			\|\rho_0^{\frac{1}{p+2}}v_0\|_{L^{p+2}}
			\le
			\|\rho_0^{\frac{1}{2}}v_0\|_{L^2}^{\frac{2}{p+2}}
			\|v_0\|_{L^\infty}^{\frac{p}{p+2}}
			\le
			\|\rho_0^{\frac{1}{2}}v_0\|_{L^2}
			+
			\|v_0\|_{L^\infty}.
		\end{aligned}
	\end{equation}
	We introduce the following notation for $k\in\mathbb N^+$:
    \[
	\begin{aligned}
	&r:=\frac{5}{3},\qquad p_k+2:=r^k,
	\qquad
	C_7:=\|\rho_0^{\frac{1}{2}}v_0\|_{L^2}+\|v_0\|_{L^\infty}+1,\\
	&\Xi(k):=\int_{Q_T}\rho|v|^{r^k}\,dxdt.
	\end{aligned}
    \]
	For every \(k\ge3\), substituting \(p=p_k\) into \eqref{x411} and using \eqref{x412} gives
	\begin{align}\label{x413}
		\Xi(k+1)
		\le
		C_6\Theta_T
		\left(
		r^{\frac{10}{3}k}\Xi(k)^r
		+
		r^{\frac{10}{3}k}
		+
		C_7^{r^{k+1}}
		\right).
	\end{align}
	Let
	\(
	\widetilde\Xi(k):=\max\{\Xi(k),C_7^{r^k}\}.
	\)
	Then assumptions \(C_6\ge1\), \(C_7\ge1\), \(\Theta_T>1\), and \(r^{\frac{10}{3}k}\ge1\) allow us to reduce \eqref{x413} to
	\begin{align}\label{x414}
		\widetilde\Xi(k+1)
		\le
		3C_6\Theta_T r^{\frac{10}{3}k}\widetilde\Xi(k)^r .
	\end{align}
	Iterating \eqref{x414} from \(\ell\) to \(k\),  for every integer \(k>\ell\), we find that
	\begin{align}\label{x415}
		\widetilde\Xi(k+1)^{\frac{1}{r^{k+1}}}
		&\le
		(3C_6\Theta_T)^{\sum_{i=\ell+1}^{k+1}r^{-i}}
		r^{\frac{10}{3}\sum_{i=\ell+1}^{k+1}(i-1)r^{-i}}
		\widetilde\Xi(\ell)^{\frac{1}{r^\ell}} .
	\end{align}
    
	We finally choose the starting index
	\[
	\ell:=\left\lfloor \log_r(\log\Theta_T)+1\right\rfloor .
	\]
	Since \(\log\Theta_T\ge\frac{25}{9}=r^2\), we have \(\ell\ge3\), and therefore
	\(p_\ell=r^\ell-2>2\). Applying Lemma~\ref{lm41} gives
	\[
	\Xi(\ell)^{\frac{1}{r^\ell}}
	\le
	(T^*)^{\frac{1}{r^\ell}}
	\sup_{0\le t\le T}
	\|\rho^{\frac{1}{r^\ell}}v\|_{L^{r^\ell}}
	\le
	(1+T^*)\widehat C r^{\frac{\ell}{2}}.
	\]
	By definition of \(\widetilde\Xi(\ell)\), we obtain
	\begin{align}\label{x416}
		\widetilde\Xi(\ell)^{\frac{1}{r^\ell}}
		\le
		(1+T^*)\widehat C C_7 r^{\frac{\ell}{2}}.
	\end{align}
	The choice of \(\ell\) further implies
	\[
	\ell \le \log_r(\log\Theta_T)+1,
	\qquad
	r^{\frac{\ell}{2}}\le r^{\frac{1}{2}}\sqrt{\log\Theta_T},
	\qquad
	\Theta_T^{\sum_{i=\ell+1}^{\infty}r^{-i}}
	=
	\Theta_T^{\frac{1}{(r-1)r^\ell}}
	\le
	e^{\frac{1}{r-1}}
	=
	e^{\frac{3}{2}}.
	\]
	Combining these estimates with \eqref{x415} and \eqref{x416}, and setting
	\[
	C_v:=(1+T^*)r^{\frac{1}{2}}\widehat C C_7(3C_6)^{\sum_{i=1}^{\infty}r^{-i}}
	r^{\frac{10}{3}\sum_{i=1}^{\infty}(i-1)r^{-i}}e^{\frac{3}{2}},
	\]
	a direct calculation shows that
	\begin{align}\label{x417}
		\left(\int_{Q_T}\rho|v|^{r^{k+1}}\,dxdt\right)^{\frac{1}{r^{k+1}}}
		\le
		C_v\sqrt{\log\Theta_T}.
	\end{align}
    It remains to remove the density weight. Since \(\rho^{-1}\le\Theta_T\), which implies that
	\[
	\left(
	\int_{Q_T}|v|^{r^{k+1}}\,dxdt
	\right)^{\frac{1}{r^{k+1}}}
	\le
	\Theta_T^{\frac{1}{r^{k+1}}}
	\left(
	\int_{Q_T}\rho|v|^{r^{k+1}}\,dxdt
	\right)^{\frac{1}{r^{k+1}}}.
	\]
	Letting \(k\to\infty\) and using \eqref{x417}, we conclude that
	\[
	\|v\|_{L_T^\infty L^\infty}
	\le
	C_v\sqrt{\log\Theta_T}.
	\]
	This completes the proof.
	\thatsall\\

    Finally, we apply the Nash--Moser iteration to establish the following uniform positive lower bound for the density.
	\begin{proposition}\label{pr43}
		There exists a constant \(C>0\), depending on
		\(T^*,q,\gamma,E_0,\Omega,\)\(\|\rho_0\|_{L^\infty},\)
		\(\|v_0\|_{L^\infty},\)
		\(\|\rho_0^{\frac{1}{q+2}}v_0\|_{L^{q+2}}\), and
		\(\|\tau_0\|_{L^\infty}\), such that
		\begin{align}\label{x418}
			\sup_{0\le t\le T}\|\tau(t)\|_{L^\infty}
			\le C.
		\end{align}
	\end{proposition}
	\pf
	The proof is divided into three steps. We first establish a short-time bound for $\Theta_T$, then extend it to $[0,2T_0]$, and finally repeat the argument finitely many times to obtain the desired estimate on $[0,T^*)$.
	
	\medskip
	\noindent
	{\bf Step 1.} We show that there exists \(T_0>0\) sufficiently small, such that
	\begin{align}\label{x419}
		\Theta_{T_0}\le C .
	\end{align}
	Recall that \(\tau=\rho^{-1}\). It follows from \eqref{x16}\(_1\) that \(\tau\) satisfies
	\begin{align}\label{x420}
		\tau_t-\Delta\tau
		+
		2\tau^{-1}|\nabla\tau|^2
		+
		v\cdot\nabla\tau
		-
		\tau\operatorname{div}v
		=0 .
	\end{align}
	For any \(p>0\), multiplying \eqref{x420} by \(\tau^{p+1}\) and integrating
	over \(\Omega\), we obtain
	\begin{equation}\label{x421}
		\begin{aligned}
			&\frac{1}{p+2}\frac{d}{dt}\int_\Omega \tau^{p+2}\,dx
			+
			(p+3)\int_\Omega \tau^p|\nabla\tau|^2\,dx
			+
			(p+3)\int_\Omega v\cdot\nabla\tau\,\tau^{p+1}\,dx\\
			=&
			\int_{\partial\Omega}\tau^{p+1}\partial_n\tau\,dS
			+
			\int_{\partial\Omega}\tau^{p+2}v\cdot n\,dS .
		\end{aligned}
	\end{equation}
	The boundary condition \(\partial_n\rho=0\) gives
	\[
	\partial_n\tau
	=
	-\rho^{-2}\partial_n\rho
	=
	0
	\quad\text{on }\partial\Omega .
	\]
	Together with \(v\cdot n=0\) from \eqref{x18}, this shows that both boundary terms in \eqref{x421} vanish. Hence
	\begin{equation}\label{x422}
		\begin{aligned}
			&\frac{1}{p+2}\frac{d}{dt}\int_\Omega \tau^{p+2}\,dx
			+
			(p+3)\int_\Omega \tau^p|\nabla\tau|^2\,dx\\
			\le&
			\frac{p+3}{2}\int_\Omega \tau^p|\nabla\tau|^2\,dx
			+
			(p+2)\|v\|_{L_T^\infty L^\infty}^2
			\int_\Omega \tau^{p+2}\,dx .
		\end{aligned}
	\end{equation}
	Integrating \eqref{x422} over \([0,T]\) and applying the Sobolev embedding together with Lemma~\ref{lm24}, where we may assume $C_{8}\ge1$ without loss of generality:
	\begin{equation}\label{x424}
		\begin{aligned}
			\int_{Q_T}\tau^{\frac{5}{3}(p+2)}\,dxdt
			&\le
			\int_0^T
			\left(\int_\Omega \tau^{p+2}\,dx\right)^{\frac{2}{3}}
			\left(\int_\Omega \tau^{3(p+2)}\,dx\right)^{\frac{1}{3}}dt\\
			&\le
			C
			\left(
			\sup_{0\le t\le T}\int_\Omega \tau^{p+2}\,dx
			+
			\int_{Q_T}\left|\nabla\tau^{\frac{p+2}{2}}\right|^2\,dxdt
			\right)^{\frac{5}{3}}\\
			&\le
			C_8(p+2)^{\frac{10}{3}}
			\|v\|_{L_T^\infty L^\infty}^{\frac{10}{3}}
			\left(\int_{Q_T}\tau^{p+2}\,dxdt\right)^{\frac{5}{3}}\\
			&\quad+
			C_8
			\left(
			\|\tau_0\|_{L^2}
			+
			\|\tau_0\|_{L^\infty}
			\right)^{\frac{5}{3}(p+2)}.
		\end{aligned}
	\end{equation}
	For \(k\ge2\), define
    \[
    \begin{aligned}
	&r:=\frac{5}{3},
	\qquad
	p_k+2:=r^k,
	\qquad
	C_9:=\|\tau_0\|_{L^2}+\|\tau_0\|_{L^\infty}+1,\\
	&\Upsilon(k):=\int_{Q_T}\tau^{r^k}\,dxdt .
	\end{aligned}
    \]
	Taking \(p=p_k\) in \eqref{x424} and using
	Proposition~\ref{pr42}, we obtain
	\begin{equation}\label{x425}
		\begin{aligned}
			\Upsilon(k+1)
			&\le
			C_8C_v^{\frac{10}{3}}
			r^{\frac{10}{3}k}
			(\log\Theta_T)^{\frac{5}{3}}
			\Upsilon(k)^r
			+
			C_8C_9^{r^{k+1}} .
		\end{aligned}
	\end{equation}
	Let
    \(
    \widetilde\Upsilon(k):=\max\{\Upsilon(k),C_9^{r^k}\}.
	\)
	As \(C_8\ge1\), \(C_9\ge1\), and \(\Theta_T>1\), estimate \eqref{x425} reduces to
	\begin{align}\label{x426}
		\widetilde\Upsilon(k+1)
		\le
		2C_8C_v^{\frac{10}{3}}
		r^{\frac{10}{3}k}
		(\log\Theta_T)^{\frac{5}{3}}
		\widetilde\Upsilon(k)^r .
	\end{align}
	Iterating \eqref{x426} from \(2\) to \(k\) yields
	\begin{equation}\label{x427}
		\begin{aligned}
			\widetilde\Upsilon(k+1)^{\frac{1}{r^{k+1}}}
			&\le
			\left(2C_8C_v^{\frac{10}{3}}\right)^{\sum_{i=1}^{\infty}r^{-i}}
			r^{\frac{10}{3}\sum_{i=1}^{\infty}(i-1)r^{-i}}
			(\log\Theta_T)^{\frac{9}{10}}
			\widetilde\Upsilon(2)^{\frac{1}{r^2}} .
		\end{aligned}
	\end{equation}
	
	It remains to estimate the initial term \(\widetilde\Upsilon(2)\). Taking
	\(
	p=r^2-2=\frac{7}{9}
	\)
	in \eqref{x422}, we get
	\begin{align}\label{x428}
		\frac{d}{dt}\int_\Omega \tau^{r^2}\,dx
		\le
		r^4\|v\|_{L_T^\infty L^\infty}^{2}
		\int_\Omega \tau^{r^2}\,dx .
	\end{align}
	By Gronwall's inequality, we obtain
	\begin{align}\label{x429}
		\sup_{0\le t\le T}\int_\Omega \tau^{r^2}\,dx
		\le
		\exp\left\{
		r^4T\|v\|_{L_T^\infty L^\infty}^{2}
		\right\}
		\int_\Omega \tau_0^{r^2}\,dx .
	\end{align}
	Combining \eqref{x429} with Proposition~\ref{pr42}, we infer that
	\begin{equation}\label{x430}
		\begin{aligned}
			\Upsilon(2)^{\frac{1}{r^2}}
			\le
			(T^*)^{\frac{1}{r^2}}
			\Theta_T^{r^2TC_v^2}
			\|\tau_0\|_{L^{r^2}} .
		\end{aligned}
	\end{equation}
	On the other hand, the assumption \(\Theta_T\ge e^{\frac{25}{9}}\) implies
	\[
	(\log\Theta_T)^{\frac{9}{10}}
	\le
	\Theta_T^{\frac{1}{2}}.
	\]
	Using this estimate together with \eqref{x427}, \eqref{x430}, and the definition of \(\widetilde\Upsilon(2)\), we obtain
	\begin{align}\label{x431}
		\left(
		\int_{Q_T}\tau^{r^{k+1}}\,dxdt
		\right)^{\frac{1}{r^{k+1}}}
		\le
		C_{10}\Theta_T^{\frac{1}{2}+r^2TC_v^2},
		\qquad k\ge2,
	\end{align}
	where 
	\[
	C_{10}:=(1+(T^*)^{\frac{1}{r^2}})C_9\left(2C_8C_v^{\frac{10}{3}}\right)^{\sum_{i=1}^{\infty}r^{-i}}
	r^{\frac{10}{3}\sum_{i=1}^{\infty}(i-1)r^{-i}}(\|\tau_0\|_{L^{r^2}}+1)
	\]
	depends only on the quantities stated in the proposition. Passing to the limit \(k\to\infty\) in \eqref{x431} and recalling the definition of \(\Theta_T\), we deduce that
	\begin{align}\label{x432}
		\Theta_T
		\le
		\left(C_{10}+e^{\frac{25}{9}}\right)
		\Theta_T^{\frac{1}{2}+r^2TC_v^2}.
	\end{align}
	Choose \(T_0>0\) sufficiently small such that
    \(
		r^2T_0C_v^2\le\frac{1}{4}.
    \)
	With \(T=T_0\), estimate \eqref{x432} becomes
	\[
	\Theta_{T_0}
	\le
	\left(C_{10}+e^{\frac{25}{9}}\right)\Theta_{T_0}^{\frac{3}{4}}
	\le
	\frac{3}{4}\Theta_{T_0}
	+
	\frac{1}{4}\left(C_{10}+e^{\frac{25}{9}}\right)^4.
	\]
	It follows that
	\begin{align}\label{x434}
		\Theta_{T_0}
		\le
		\left(C_{10}+e^{\frac{25}{9}}\right)^4 .
	\end{align}
    This proves \eqref{x419}.
	
	\medskip
	\noindent
	{\bf Step 2.} We extend the bound obtained in Step~1 to the interval \([0,2T_0]\). We may assume that \(2T_0<T^*\). Fix any \(T\in(T_0,2T_0)\). Applying \eqref{x428} on \([T_0,T]\), we obtain
	\begin{align}\label{x435}
		\sup_{T_0\le t\le T}\int_\Omega \tau^{r^2}\,dx
		\le
		\exp\left\{
		r^4T_0\|v\|_{L_T^\infty L^\infty}^{2}
		\right\}
		\int_\Omega \tau^{r^2}(T_0)\,dx .
	\end{align}
	Combining \eqref{x435} with Proposition~\ref{pr42} gives
	\begin{equation}\label{x436}
		\begin{aligned}
			\Upsilon(2)^{\frac{1}{r^2}}
			\le
			(T^*)^{\frac{1}{r^2}}\Theta_{T_0}
			+
			(T^*)^{\frac{1}{r^2}}
			\Theta_T^{\frac{1}{4}}\|\tau(T_0)\|_{L^{r^2}} .
		\end{aligned}
	\end{equation}
	Substituting \eqref{x436} into \eqref{x427}, we find that, for every
	\(k\ge2\),
	\begin{equation}\label{x437}
		\begin{aligned}
			\left(
			\int_{Q_T}\tau^{r^{k+1}}\,dxdt
			\right)^{\frac{1}{r^{k+1}}}
			&\le
			\left(2C_8C_v^{\frac{10}{3}}\right)^{\sum_{i=1}^{\infty}r^{-i}}
			r^{\frac{10}{3}\sum_{i=1}^{\infty}(i-1)r^{-i}}
			(\log\Theta_T)^{\frac{9}{10}}\\
			&\quad\times
			\max\left\{
			(T^*)^{\frac{1}{r^2}}\Theta_{T_0}
			+
			(T^*)^{\frac{1}{r^2}}
			\Theta_T^{\frac{1}{4}}\|\tau(T_0)\|_{L^{r^2}},
			C_9
			\right\}\\
			&\le
			C_{11}\Theta_T^{\frac{3}{4}},
		\end{aligned}
	\end{equation}
	where
	\[
	C_{11}:=2\Theta_{T_0}C_{10}
	\]
	depends only on the quantities stated in the proposition. Passing to the limit \(k\to\infty\) in \eqref{x437} yields
	\[
	\|\tau\|_{L_T^\infty L^\infty}
	\le
	C_{11}\Theta_T^{\frac{3}{4}}.
	\]
	Recalling the definition of \(\Theta_T\), we therefore obtain
	\[
	\Theta_T
	\le
	\left(C_{11}+e^{\frac{25}{9}}\right)\Theta_T^{\frac{3}{4}}.
	\]
	By Young's inequality and the arbitrariness of \(T\in(T_0,2T_0)\), we conclude that
	\begin{align}\label{x438}
		\Theta_{2T_0}
		\le
		\left(C_{11}+e^{\frac{25}{9}}\right)^4 .
	\end{align}
	
	\medskip
	\noindent
	{\bf Step 3.} Choose \(j\in\mathbb N^+\) such that 
	\[
	jT_0\le T^*<(j+1)T_0.
	\]
	Repeating Step 2 on the finite sequence of intervals
	\(
	[0,T_0],\ [T_0,2T_0],\ \cdots,\ [jT_0,T^*)
	\)
	yields
	\[
	\Theta_T\le C,
	\qquad
	\forall\,0<T<T^*.
	\]
	Recalling the definition of \(\Theta_T\), we obtain
	\[
	\sup_{0\le t\le T}\|\tau(t)\|_{L^\infty}
	\le C,
	\qquad
	\forall\,0<T<T^*.
	\]
	This proves \eqref{x418}, which completes the proof.
	\thatsall
	
	\section{High-order estimates (\(N=2,3\))}
	
	With the upper and lower bounds of the density at hand, we are ready to derive higher-order estimates for \((\rho,v)\). Unlike the periodic case, one has to use the Neumann elliptic estimate for the density and the elliptic estimate for the effective velocity under the slip boundary condition. These estimates were prepared in Section~2 and will be used repeatedly.
	Since every such embedding used below remains valid in dimension two (2D), the same proof covers \(N=2\).

    We begin with the following $H^2$-estimate for $\rho$ and $H^1$-estimate for $v$, where Lemmas~\ref{lm27} and~\ref{lm28} are used to control $\|\rho\|_{H^3}$ and $\|v\|_{H^2}$.
	
	\begin{proposition}\label{pr51}
		There exists a constant \(C>0\), depending on
		\(T^*,\gamma,\Omega,\)\(\|\rho_0\|_{L^\infty},\)
		\(\|v_0\|_{L^\infty},\)
		\(\|\tau_0\|_{L^\infty},\)
		\(\|\rho_0\|_{H^{2}}\) and \(\|v_0\|_{H^{1}}\)
		, such that
		\begin{align}\label{x51}
			\sup_{0\le t\le T}
			\left(
			\|\rho\|_{H^2}
			+
			\|\rho_t\|_{L^2}
			+
			\|v\|_{H^1}
			\right)
			+
			\int_0^T
			\left(
			\|\rho\|_{H^3}^2
			+
			\|\rho_t\|_{H^1}^2
			+
			\|v\|_{H^2}^2
			+
			\|v_t\|_{L^2}^2
			\right)\,dt
			\le C.
		\end{align}
	\end{proposition}
	\pf
	Since \(\rho\) has a positive lower bound, system \eqref{x16} can be written as
	\begin{align}\label{x52}
		\begin{cases}
			\rho_t+\operatorname{div}(\rho v)-\Delta\rho=0,\\[2mm]
			v_t+(v-2\nabla\log\rho)\cdot\nabla v-\Delta v
			+\nabla\Pi_\gamma(\rho)=0.
		\end{cases}
	\end{align}
	
	We begin with the estimate for the density. Multiplying \eqref{x52}\(_1\) by \(-\Delta\rho_t\), integrating over \(\Omega\), and using \eqref{x18}, we obtain
	\begin{equation}\label{x53}
		\begin{aligned}
			\frac{1}{2}\frac{d}{dt}\int_\Omega |\Delta\rho|^2\,dx
			+
			\int_\Omega |\nabla\rho_t|^2\,dx
			&=
			-\int_\Omega \nabla\operatorname{div}(\rho v)\cdot\nabla\rho_t\,dx\\
			&\le
			\frac{1}{4}\int_\Omega |\nabla\rho_t|^2\,dx
			+
			C\int_\Omega |\nabla\operatorname{div}(\rho v)|^2\,dx .
		\end{aligned}
	\end{equation}
	Applying \(\nabla\) to \eqref{x52}\(_1\) and estimating the resulting terms gives 
	\begin{align}\label{x54}
		\frac{1}{8}\int_\Omega |\nabla\Delta\rho|^2\,dx
		\le
		\frac{1}{4}\int_\Omega |\nabla\rho_t|^2\,dx
		+
		\frac{1}{4}\int_\Omega |\nabla\operatorname{div}(\rho v)|^2\,dx .
	\end{align}
	Combining \eqref{x53} and \eqref{x54}, we have
	\begin{align}\label{x55}
		\frac{1}{2}\frac{d}{dt}\int_\Omega |\Delta\rho|^2\,dx
		+
		\frac{1}{2}\int_\Omega |\nabla\rho_t|^2\,dx
		+
		\frac{1}{8}\int_\Omega |\nabla\Delta\rho|^2\,dx
		\le
		C\int_\Omega |\nabla\operatorname{div}(\rho v)|^2\,dx .
	\end{align}
    
	It remains to estimate the term on the right-hand side of \eqref{x55}. Using Lemma~\ref{lm22}, \ref{lm27}, \ref{lm28}, Proposition~\ref{pr34}, \ref{pr42} and \ref{pr43}, we obtain
	\begin{equation}\label{x56}
		\begin{aligned}
			\int_\Omega |\nabla\operatorname{div}(\rho v)|^2\,dx
			&\le
			C\int_\Omega |\nabla^2\rho|^2|v|^2\,dx
			+
			C\int_\Omega |\nabla\rho|^2|\nabla v|^2\,dx\\
			&\quad+
			C\int_\Omega \rho^2|\nabla^2v|^2\,dx\\
			&\le
			C\int_\Omega |\Delta\rho|^2\,dx
			+
			C\int_\Omega |\nabla\rho|^2|\nabla v|^2\,dx
			+
			C_8\int_\Omega |\Delta v|^2\,dx\\
			&\quad+
			C\int_\Omega |\rho|^2\,dx
			+
			C\int_\Omega |v|^2\,dx\\
			&\le
			C\int_\Omega |\Delta\rho|^2\,dx
			+
			C\int_\Omega |\nabla\rho|^2|\nabla v|^2\,dx
			+
			C_8\int_\Omega |\Delta v|^2\,dx
			+
			C,
		\end{aligned}
	\end{equation}
	where in the last inequality, we have used the boundedness of \(\rho\) and \(v\), as well as the boundedness of \(\Omega\). Substituting \eqref{x56} into \eqref{x55} yields
	\begin{equation}\label{x57}
		\begin{aligned}
			&\frac{1}{2}\frac{d}{dt}\int_\Omega |\Delta\rho|^2\,dx
			+
			\frac{1}{2}\int_\Omega |\nabla\rho_t|^2\,dx
			+
			\frac{1}{8}\int_\Omega |\nabla\Delta\rho|^2\,dx\\
			\le&
			C\int_\Omega |\Delta\rho|^2\,dx
			+
			C\int_\Omega |\nabla\rho|^2|\nabla v|^2\,dx
			+
			C_8\int_\Omega |\Delta v|^2\,dx
			+
			C .
		\end{aligned}
	\end{equation}
	
	We now turn to the estimate for \(v\). The boundary contribution
	arising from the diffusion term is treated as in
	Proposition~\ref{pr32}. In fact, \eqref{x3111} and the boundary conditions \eqref{x18}
	gives
	\[
	-\int_\Omega v_t\cdot\Delta v\,dx
	=
	\frac{1}{2}\frac{d}{dt}
	\left(
	\int_\Omega |\nabla v|^2\,dx
	+
	\int_{\partial\Omega}v\cdot\kappa(x)v\,dS
	\right).
	\]
	Testing \eqref{x52}\(_2\) by \(-\Delta v\) gives
	\begin{equation}\label{x58}
		\begin{aligned}
			&\frac{1}{2}\frac{d}{dt}
			\left(
			\int_\Omega |\nabla v|^2\,dx
			+
			\int_{\partial\Omega}v\cdot\kappa(x)v\,dS
			\right)
			+
			\int_\Omega |\Delta v|^2\,dx\\
			&\qquad\le
			\int_\Omega |\Delta v|
			\left(
			|v||\nabla v|
			+
			|\nabla\log\rho||\nabla v|
			+
			|\nabla\rho^{\gamma-1}|
			\right)\,dx\\
			&\qquad\le
			\frac{1}{2}\int_\Omega |\Delta v|^2\,dx
			+
			C\int_\Omega |\nabla v|^2\,dx
			+
			C\int_\Omega |\nabla\rho|^2|\nabla v|^2\,dx
			+
			C .
		\end{aligned}
	\end{equation}
	Multiplying \eqref{x58} by \(4C_8\) and adding it to
	\eqref{x57}, we obtain
	\begin{equation}\label{x59}
		\begin{aligned}
			&\frac{1}{2}\frac{d}{dt}\int_\Omega |\Delta\rho|^2\,dx
			+
			2C_8\frac{d}{dt}
			\left(
			\int_\Omega |\nabla v|^2\,dx
			+
			\int_{\partial\Omega}v\cdot\kappa(x)v\,dS
			\right)\\
			&+
			\frac{1}{2}\int_\Omega |\nabla\rho_t|^2\,dx
			+
			\frac{1}{8}\int_\Omega |\nabla\Delta\rho|^2\,dx
			+
			C_8\int_\Omega |\Delta v|^2\,dx\\
			\le&
			C\int_\Omega |\Delta\rho|^2\,dx
			+
			C\int_\Omega |\nabla v|^2\,dx
			+
			C\int_\Omega |\nabla\rho|^2|\nabla v|^2\,dx
			+
			C .
		\end{aligned}
	\end{equation}
	
	The only term that still requires additional control is
	\begin{equation}\label{x5101}
	\int_\Omega |\nabla\rho|^2|\nabla v|^2\,dx .
	\end{equation}
    Integration by parts yields
	\begin{equation}\label{x510}
		\begin{aligned}
			\int_\Omega |\nabla\rho|^2|\nabla v|^2\,dx
			&=
			\int_\Omega \partial_i\rho\,\partial_i\rho\,\partial_jv^k\,\partial_jv^k\,dx\\
			&=
			-2\int_\Omega \partial_{ij}\rho\,\partial_i\rho\,v^k\partial_jv^k\,dx
			-
			\int_\Omega |\nabla\rho|^2v^k\Delta v^k\,dx\\
			&\quad
			+
			\int_{\partial\Omega}|\nabla\rho|^2(n\cdot\nabla v)\cdot v\,dS .
		\end{aligned}
	\end{equation}
	The assumption \(\kappa(x)\ge0\) therefore shows that the boundary contribution in \eqref{x510} is negative. Lemmas~\ref{lm22}, \ref{lm24}, \ref{lm27}, and \ref{lm28}, together with Proposition~\ref{pr42}, H\"older's inequality, and Young's inequality, then imply that, for any \(\varepsilon_1,\varepsilon_2>0\),
	\begin{equation}\label{x511}
	\begin{aligned}
		\int_\Omega |\nabla\rho|^2|\nabla v|^2\,dx
		&\le
		C\|\nabla^2\rho\|_{L^6}\|\nabla\rho\|_{L^2}
		\|v\|_{L^\infty}\|\nabla v\|_{L^3}
		+
		C\|\nabla\rho\|_{L^4}^2
		\|v\|_{L^\infty}\|\Delta v\|_{L^2}\\
		&\le
		\varepsilon_1\|\nabla^2\rho\|_{H^1}^2
		+
		\varepsilon_2\|\nabla v\|_{H^1}^2
		+
		C(\varepsilon_1,\varepsilon_2)\|\nabla v\|_{L^2}^2\\
		&\quad+
		C(\varepsilon_1,\varepsilon_2)\|\Delta\rho\|_{L^2}^2
		+
		C(\varepsilon_1,\varepsilon_2).
	\end{aligned}
    \end{equation}
    By Lemmas~\ref{lm27} and~\ref{lm28}, we may choose $\varepsilon_1,\varepsilon_2>0$ small enough, depending only on the corresponding elliptic constants. Substituting \eqref{x511} into \eqref{x59}, a direct calculation shows that
	\begin{equation}\label{x512}
		\begin{aligned}
			&\frac{1}{2}\frac{d}{dt}\int_\Omega |\Delta\rho|^2\,dx
			+
			2C_8\frac{d}{dt}
			\left(
			\int_\Omega |\nabla v|^2\,dx
			+
			\int_{\partial\Omega}v\cdot\kappa(x)v\,dS
			\right)\\
			&+
			\frac{1}{2}\int_\Omega |\nabla\rho_t|^2\,dx
			+
			\frac{1}{16}\int_\Omega |\nabla\Delta\rho|^2\,dx
			+
			\frac{C_8}{2}\int_\Omega |\Delta v|^2\,dx\\
			\le&
			C\int_\Omega |\Delta\rho|^2\,dx
			+
			C\int_\Omega |\nabla v|^2\,dx
			+
			C .
		\end{aligned}
	\end{equation}
    
	The boundary part of the energy in \eqref{x512} is positive, since \(\kappa(x)\ge0\). The trace theorem and the regularity of the initial data also give
	\[
	\int_{\partial\Omega}v_0\cdot\kappa(x)v_0\,dS
	\le
	C\|v_0\|_{H^1}^2.
	\]
	Applying Gronwall's inequality to \eqref{x512} yields
	\begin{align}\label{x513}
		\sup_{0\le t\le T}
		\left(
		\|\Delta\rho\|_{L^2}^2
		+
		\|\nabla v\|_{L^2}^2
		\right)
		+
		\int_0^T
		\left(
		\|\nabla\rho_t\|_{L^2}^2
		+
		\|\nabla\Delta\rho\|_{L^2}^2
		+
		\|\Delta v\|_{L^2}^2
		\right)\,dt
		\le C .
	\end{align}
	By Lemma~\ref{lm27} and \ref{lm28}, it follows that
	\begin{align}\label{x514}
		\sup_{0\le t\le T}
		\left(
		\|\rho\|_{H^2}^2
		+
		\|v\|_{H^1}^2
		\right)
		+
		\int_0^T
		\left(
		\|\rho\|_{H^3}^2
		+
		\|\rho_t\|_{H^1}^2
		+
		\|v\|_{H^2}^2
		\right)\,dt
		\le C .
	\end{align}
	
	To complete the estimate, it remains to control the time derivatives. From \eqref{x52}\(_1\), Lemma~\ref{lm22}, Propositions~\ref{pr34} and \ref{pr42}, and \eqref{x514}, we have
	\begin{equation}\label{x515}
		\begin{aligned}
			\|\rho_t\|_{L^2}
			&\le
			C\left(\|\Delta\rho\|_{L^2}
			+
			\|\nabla(\rho v)\|_{L^2}\right)\\
			&\le
			C\left(\|\Delta\rho\|_{L^2}
			+
			\|\rho\|_{L^\infty}\|\nabla v\|_{L^2}
			+
			\|\nabla\rho\|_{L^2}\|v\|_{L^\infty}\right)
			\le C .
		\end{aligned}
	\end{equation}
	Similarly, \eqref{x52}\(_2\), Lemmas~\ref{lm22}, \ref{lm24}, and \ref{lm28}, together with \eqref{x514} gives
	\begin{equation}\label{x516}
		\begin{aligned}
			\|v_t\|_{L^2}
			&\le
			C\left(\||v||\nabla v|\|_{L^2}
			+
			\||\nabla\log\rho||\nabla v|\|_{L^2}
			+
			\|\Delta v\|_{L^2}
			+
			\|\nabla\rho^{\gamma-1}\|_{L^2}\right)\\
			&\le
			C\left(\|v\|_{L^\infty}\|\nabla v\|_{L^2}
			+
			\|\nabla\rho\|_{L^6}\|\nabla v\|_{L^3}
			+
			\|\Delta v\|_{L^2}
			+
			\|\nabla\rho\|_{L^2}\right)\\
			&\le
			C+C\|\nabla^2 v\|_{L^2}.
		\end{aligned}
	\end{equation}
	Combining \eqref{x514}, \eqref{x515}, and \eqref{x516} together, we obtain \eqref{x51}, which completes the proof.
	\thatsall\\

    Next, we establish the following $H^2$-estimate for $v$ and $L^2$-estimate for $v_t$, using the elliptic estimate in Lemma~\ref{lm28} to control $\|v\|_{H^3}$.
	\begin{proposition}\label{pr52}
		There exists a constant \(C>0\), depending on
		\(T^*,\gamma,\Omega,\|\tau_0\|_{L^\infty},
		\|\rho_0\|_{H^2}\) and \(\|v_0\|_{H^2}\), such that
		\begin{align}\label{x517}
			\sup_{0\le t\le T}
			\left(
			\|v\|_{H^2}
			+
			\|v_t\|_{L^2}
			\right)
			+
			\int_0^T
			\left(
			\|v\|_{H^3}^2
			+
			\|v_t\|_{H^1}^2
			\right)\,dt
			\le C .
		\end{align}
	\end{proposition}
	\pf
	Differentiating \eqref{x52}\(_2\) with respect to \(t\), we obtain
	\begin{align}\label{x518}
		v_{tt}
		+
		(v-2\nabla\log\rho)_t\cdot\nabla v
		+
		(v-2\nabla\log\rho)\cdot\nabla v_t
		-
		\Delta v_t
		+
		\frac{\gamma}{\gamma-1}
		\bigl(\nabla\rho^{\gamma-1}\bigr)_t
		=0 .
	\end{align}
	We first examine the boundary terms arising in the corresponding \(L^2\)-estimate. Since
    \[
    v_t\cdot n=0,\qquad
    \operatorname{curl}v_t\times n=0
    \quad\text{on }\partial\Omega,
    \]
    the same boundary identity as in \eqref{x3121} gives
	\[
	-\int_\Omega \Delta v_t\cdot v_t\,dx
	=
	\int_\Omega |\nabla v_t|^2\,dx
	-
	\int_{\partial\Omega}(n\cdot\nabla v_t)\cdot v_t\,dS
	\ge
	\int_\Omega |\nabla v_t|^2\,dx .
	\]
	After integrating by parts, the remaining boundary contribution is
	\[
	\mathcal B_t
	=
	2\int_{\partial\Omega}(\log\rho)_t(n\cdot\nabla v)\cdot v_t\,dS .
	\]
	Since \(v_t\cdot n=0\) on \(\partial\Omega\), the vector \(v_t\) is tangential there. Hence \eqref{x3111} yields
	\[
	(n\cdot\nabla v)\cdot v_t
	=
	(v_t\cdot\nabla v)\cdot n
	=
	-(v_t\cdot\nabla n)\cdot v.
	\]
	Using above equation, the trace theorem, H\"older's inequality and Young's inequality, we have, for any
	\(\varepsilon>0\),
	\begin{equation}\label{x519}
		\begin{aligned}
			|\mathcal B_t|
			&\le
			C\|(\log\rho)_t\|_{L^2(\partial\Omega)}
			\|v\|_{L^\infty(\partial\Omega)}
			\|v_t\|_{L^2(\partial\Omega)} \\
			&\le
			\varepsilon\|v_t\|_{H^1}^2
			+
			C(\varepsilon)\|\rho_t\|_{H^1}^2
			+
			C(\varepsilon)\|v_t\|_{L^2}^2 .
		\end{aligned}
	\end{equation}
    Testing \eqref{x518} by \(v_t\), integrating by parts, and using \eqref{x519} with \(\varepsilon=\frac18\), together with Young's inequality, we obtain
	\begin{equation}\label{x520}
		\begin{aligned}
			&\frac{1}{2}\frac{d}{dt}\int_\Omega |v_t|^2\,dx
			+
			\int_\Omega |\nabla v_t|^2\,dx \\
			\le&
			\frac{1}{4}\int_\Omega |\nabla v_t|^2\,dx
			+
			C\int_\Omega
			\left(
			|v|^2|v_t|^2
			+
			|\rho_t|^2|\nabla v|^2
			+
			|\nabla\rho|^2|v_t|^2
			+
			|\rho_t|^2
			\right)\,dx \\
			&+
			C\int_\Omega |\rho_t||\Delta v||v_t|\,dx
			+
			C\left(\|\rho_t\|_{H^1}^2
			+
			\|v_t\|_{L^2}^2\right) .
		\end{aligned}
	\end{equation}
	
	We next recover the third-order spatial derivatives of \(v\) from \eqref{x52}\(_2\). Differentiating the equation in space gives
	\[
	\nabla\Delta v
	=
	\nabla v_t
	+
	\nabla\bigl((v-2\nabla\log\rho)\cdot\nabla v\bigr)
	+
	\frac{\gamma}{\gamma-1}
	\nabla^2\rho^{\gamma-1}.
	\]
	Hence, using the upper and lower bounds of \(\rho\), Lemma~\ref{lm24}, \ref{lm22} and \ref{lm27}, we have
	\begin{equation}\label{x521}
		\begin{aligned}
			\frac{1}{12}\int_\Omega |\nabla\Delta v|^2\,dx
			&\le
			\frac{1}{4}\int_\Omega |\nabla v_t|^2\,dx\\
			&\quad
			+
			C\int_\Omega
			\left(
			|\nabla v|^4
			+
			|v|^2|\nabla^2v|^2
			+
			|\nabla\rho|^4|\nabla v|^2
			+
			|\nabla^2\rho|^2|\nabla v|^2
			\right)\,dx\\
			&\quad
			+
			C\int_\Omega
			\left(
			|\nabla\rho|^2|\nabla^2v|^2
			+
			|\nabla\rho|^4
			+
			|\nabla^2\rho|^2
			\right)\,dx .
		\end{aligned}
	\end{equation}
	Adding \eqref{x520} and \eqref{x521}, and applying Lemmas~\ref{lm22}, \ref{lm24}, \ref{lm27}, and \ref{lm28}, Proposition~\ref{pr51}, and Young's inequality, a direct calculation shows that
	\begin{equation}\label{x522}
		\begin{aligned}
			&\frac{1}{2}\frac{d}{dt}\int_\Omega |v_t|^2\,dx
			+
			\frac{1}{2}\int_\Omega |\nabla v_t|^2\,dx
			+
			\frac{1}{12}\int_\Omega |\nabla\Delta v|^2\,dx \\
			\le&
			C\left(
			\|v_t\|_{L^2}^2
			+
			\|\nabla v\|_{L^2}^{\frac{1}{2}}\|\nabla v\|_{H^2}^{\frac{3}{2}}
			+
			\|\nabla\rho\|_{H^2}^2\|v_t\|_{L^2}^2
			+
			1
			\right) \\
			&
			+
			C\|\nabla^2v\|_{L^2}^{\frac{1}{2}}\|\nabla^2v\|_{H^1}^{\frac{1}{2}}
			\|v_t\|_{H^1}
			+
			C\|\nabla v\|_{L^2}^{\frac{3}{2}}\|\nabla v\|_{H^1}^{\frac{1}{2}}
			\|\nabla^2v\|_{H^1} \\
			&
			+
			C\|\nabla^2v\|_{L^2}^2
			+
			C\|\nabla v\|_{L^2}^{\frac{1}{2}}\|\nabla v\|_{H^2}^{\frac{3}{2}}
			+
			C\|\nabla^2v\|_{L^2}\|\nabla^2v\|_{H^1}
			+
			C\|\rho_t\|_{H^1}^2
			+
			C\\
			\le&
			\frac{1}{4}\|v_t\|_{H^1}^2
			+
			\varepsilon\|\nabla v\|_{H^2}^2
			+
			C(\varepsilon)\left(
			1
			+
			\|\nabla\rho\|_{H^2}^2
			\right)\|v_t\|_{L^2}^2\\
			&+
			C(\varepsilon)\left(
			1
			+
			\|\nabla v\|_{L^2}^2
			+
			\|\nabla^2v\|_{L^2}^2
			+
			\|\rho_t\|_{H^1}^2
			\right).
		\end{aligned}
	\end{equation}
    Here \(\varepsilon>0\) is arbitrary. In obtaining the last line of \eqref{x522}, the linear high-order term is estimated once more by Young's inequality:
	\[
	C\|\nabla^2v\|_{L^2}\|\nabla^2v\|_{H^1}
	\le
	\varepsilon\|\nabla v\|_{H^2}^2+C_{\varepsilon}\|\nabla^2v\|_{L^2}^2,
	\]
    To absorb the remaining high-order terms, we use Lemma~\ref{lm28} with \(k=1\). There exists an constant \(C_0>0\) such that
	\[
	\|\nabla v\|_{H^2}^2
	\le
	C_0\left(
	\|\nabla\Delta v\|_{L^2}^2
	+
	\|\Delta v\|_{L^2}^2
	+
	\|v\|_{L^2}^2
	\right).
	\]
	Choose \(\varepsilon>0\) small enough so that \(\varepsilon C_0<\frac{1}{24}\), we have
	\begin{equation}\label{x523}
		\begin{aligned}
			&\frac{1}{2}\frac{d}{dt}\int_\Omega |v_t|^2\,dx
			+
			\frac{1}{4}\int_\Omega |\nabla v_t|^2\,dx
			+
			\frac{1}{24}\int_\Omega |\nabla\Delta v|^2\,dx\\
			\le&
			C\left(1+\|\nabla\rho\|_{H^2}^2\right)\|v_t\|_{L^2}^2
			+
			C\left(
			1+\|\nabla v\|_{L^2}^2+\|\nabla^2v\|_{L^2}^2
			+
			\|\rho_t\|_{H^1}^2
			\right).
		\end{aligned}
	\end{equation}
	By Proposition~\ref{pr51}, the coefficients on the right-hand side are integrable in time. Hence Gronwall's inequality gives
	\begin{align}\label{x524}
		\sup_{0\le t\le T}\|v_t\|_{L^2}^2
		+
		\int_0^T
		\left(
		\|v_t\|_{H^1}^2
		+
		\|v\|_{H^3}^2
		\right)\,dt
		\le C,
	\end{align}
	where we use Lemma~\ref{lm28} to pass from
	\(\|\nabla\Delta v\|_{L^2}\) to \(\|v\|_{H^3}\).
	
	It remains to establish the uniform-in-time \(H^2\)-bound for \(v\). Using \eqref{x52}\(_2\), Propositions~\ref{pr42} and \ref{pr51}, together with Lemmas~\ref{lm22}, \ref{lm24} and \ref{lm28}, we have
	\begin{equation}\label{x525}
		\begin{aligned}
			\|\nabla^2v\|_{L^2}
			&\le
			C\|\Delta v\|_{L^2}
			+
			C\|v\|_{L^2} \\
			&\le
			C\left(
			\|v_t\|_{L^2}
			+
			\|v\|_{L^\infty}\|\nabla v\|_{L^2}
			+
			\|\nabla\log\rho\|_{L^6}\|\nabla v\|_{L^3}
			+
			\|\nabla\rho\|_{L^2}
			\right)
			+
			C\\
			&\le
			\frac{1}{2}\|\nabla^2v\|_{L^2}
			+
			C.
		\end{aligned}
	\end{equation}
	Absorbing the first term on the right-hand side, we obtain
	\begin{align}\label{x526}
		\sup_{0\le t\le T}\|v\|_{H^2}\le C .
	\end{align}
	Combining \eqref{x524} and \eqref{x526} proves \eqref{x517} and completes the proof.
	\thatsall\\

    Finally, we derive the following $H^3$-estimate for $\rho$ and $H^1$-estimate for $\rho_t$, using the elliptic estimate in Lemma~\ref{lm27} to control $\|\rho\|_{H^4}$ under the Neumann boundary condition \eqref{x15}.
	\begin{proposition}\label{pr53}
		There exists a constant \(C>0\), depending on
		\(T^*,\gamma,\Omega,\|\tau_0\|_{L^\infty},
		\|\rho_0\|_{H^3}\) and \(\|v_0\|_{H^2}\), such that
		\begin{align}\label{x527}
			&\sup_{0\le t\le T}
			\left(
			\|\rho\|_{H^3}
			+
			\|\rho_t\|_{H^1}
			\right)
			+
			\int_0^T
			\left(
			\|\rho\|_{H^4}^2
			+
			\|\rho_t\|_{H^2}^2
			+
			\|\rho_{tt}\|_{L^2}^2
			\right)\,dt
			\le C.
		\end{align}
	\end{proposition}
	\pf
	Differentiating \eqref{x52}\(_1\) with respect to $t$ gives
	\begin{align}\label{x528}
		\rho_{tt}-\Delta \rho_t=-\operatorname{div}(\rho_tv+\rho v_t),
		\qquad \partial_n \rho_t=0\quad\text{on }\partial\Omega.
	\end{align}
	By the Sobolev embeddings for \(N=2,3\), Proposition~\ref{pr51}, \ref{pr52}, we have
	\begin{equation}\label{x529}
		\begin{aligned}
			\|\operatorname{div}(\rho_tv+\rho v_t)\|_{L^2}^2
			&\le C\|v\|_{L^\infty}^2\|\nabla \rho_t\|_{L^2}^2
			+C\|\rho_t\|_{L^2}^2\|\nabla v\|_{L^\infty}^2 \\
			&\quad+C\|\rho\|_{L^\infty}^2\|\nabla v_t\|_{L^2}^2
			+C\|\nabla\rho\|_{L^\infty}^2\|v_t\|_{L^2}^2 \\
			&\le C\|\nabla \rho_t\|_{L^2}^2+F(t).
		\end{aligned}
	\end{equation}
	where, in view of Propositions~\ref{pr51} and~\ref{pr52}, we may take
	\[
	F(t)=C\left(\|v\|_{H^3}^2+\|v_t\|_{H^1}^2+\|\rho\|_{H^3}^2\right)\in L^1(0,T).
	\]
	Multiplying \eqref{x528} by \(\rho_{tt}\), integrating over \(\Omega\), and applying
	Young's inequality yields
	\[
	\frac{d}{dt}\|\nabla \rho_t\|_{L^2}^2+\|\rho_{tt}\|_{L^2}^2
	\le C\|\nabla \rho_t\|_{L^2}^2+F(t).
	\]
	The initial value satisfies
	\[
	\rho_t(0)=\Delta\rho_0-\operatorname{div}(\rho_0v_0)\in H^1(\Omega),
	\]
	Gronwall's inequality therefore gives the uniform \(H^1\)-bound for \(\rho_t\), together with
    \(
    \rho_{tt}\in L^2(0,T;L^2).
    \)
    Returning to \eqref{x528} and using \eqref{x529}, we further obtain
    \[
    \Delta \rho_t\in L^2(0,T;L^2).
    \]
    Employing Lemma~\ref{lm27} yields
    \(
    \rho_t\in L^2(0,T;H^2).
    \)
    Hence
	\begin{align}\label{x530}
		\sup_{0\le t\le T}\|\rho_t\|_{H^1}^2
		+
		\int_0^T
		\left(
		\|\rho_{tt}\|_{L^2}^2
		+\|\rho_t\|_{H^2}^2
		\right)\,dt
		\le C.
	\end{align}
	
	Next, we recover the higher spatial regularity of \(\rho\). Since the transport equation can be rewritten as
	\begin{align}\label{x531}
		\Delta\rho=\rho_t+\operatorname{div}(\rho v),
	\end{align}
	and the standard Sobolev product estimate gives
	\begin{align}\label{x532}
		\|\operatorname{div}(\rho v)\|_{H^2}
		\le C\left(
		\|\rho\|_{H^3}\|v\|_{H^2}
		+\|\rho\|_{H^2}\|v\|_{H^3}
		\right),
	\end{align}
	Propositions~\ref{pr51}, \ref{pr52}, \eqref{x530}, \eqref{x532}, together with Lemma~\ref{lm27} yield
	\begin{align}\label{x533}
		\int_0^T\|\rho\|_{H^4}^2\,dt\le C.
	\end{align}
	
    To obtain the corresponding bound on \([0,T]\), we use the algebra property of \(H^2(\Omega)\) for \(N=2,3\):
	\[
	\|\operatorname{div}(\rho v)\|_{H^1}
	\le
	C\|\rho\|_{H^2}\|v\|_{H^2}
	\le C.
	\]
	Together with \eqref{x530} and \eqref{x531}, Lemma~\ref{lm27} with
	\(k=1\) gives
	\[
	\sup_{0\le t\le T}\|\rho\|_{H^3}
	\le C,
	\]
	which combined with the preceding estimates proves \eqref{x527}. This completes the proof.
	\thatsall
	
	\section{Proof of Theorem~\ref{th11}}
	
	Let \(T^*\) be the maximal existence time of the local strong solution. We first consider the 3D case. By the a priori estimates obtained in the previous sections, for any \(0<T<T^*\), we have
	\begin{equation}\label{x61}
		\begin{aligned}
			&\rho\in L^\infty(0,T;H^3)\cap L^2(0,T;H^4),
			\qquad
			\rho_t\in L^\infty(0,T;H^1)\cap L^2(0,T;H^2),\\
			&v\in L^\infty(0,T;H^2)\cap L^2(0,T;H^3),
			\qquad
			v_t\in L^\infty(0,T;L^2)\cap L^2(0,T;H^1).
		\end{aligned}
	\end{equation}
	Propositions~\ref{pr34} and~\ref{pr43} also gives
	\begin{align}\label{x62}
		0<c\le \rho(x,t)\le C,
		\qquad
		(x,t)\in \Omega\times[0,T],
	\end{align}
	where the constants are independent of \(T<T^*\).
	
	We next recover the time continuity required in the definition of a strong solution.By \cite[Section 5.9.2, Theorem 4]{EvansPDE}, we use the following standard consequence
	\[
	f\in L^2(0,T;H^{s+1}),
	\qquad
	f_t\in L^2(0,T;H^{s-1})
	\quad\Longrightarrow\quad
	f\in C([0,T];H^s).
	\]
	Applying this result to \(\rho\) with \(s=3\), to \(v\) with \(s=2\), and to \(\rho_t\) with \(s=1\), we obtain
	\begin{align}\label{x63}
		\rho\in C([0,T];H^3),
		\qquad
		v\in C([0,T];H^2),
		\qquad
		\rho_t\in C([0,T];H^1).
	\end{align}
	Since \(H^3(\Omega)\hookrightarrow C^1(\overline\Omega)\) for \(N=3\), and
	\(\rho\) has a positive lower bound, the composition property of Sobolev spaces gives
	\[
	\nabla\log\rho\in C([0,T];H^2)\cap L^2(0,T;H^3).
	\]
	Recalling that
	\(
	v=u+\nabla\log\rho,
	\)
	we therefore have
	\begin{align}\label{x64}
		u=v-\nabla\log\rho
		\in C([0,T];H^2).
	\end{align}
	For the time derivative of \(u\), we write
	\[
	u_t
	=
	v_t-\nabla\left(\frac{\rho_t}{\rho}\right),
	\]
	Using \eqref{x61}, \eqref{x62}, and the Sobolev product estimates in \(H^2(\Omega)\), we further obtain
	\begin{align}\label{x65}
		u_t\in L^\infty(0,T;L^2)\cap L^2(0,T;H^1).
	\end{align}
	
	We are now in a position to close the continuation argument. Suppose, to the contrary, that \(T^*<\infty\). Since the constants in the above estimates can be bounded in terms of \(T^*\), all the preceding bounds hold uniformly for \(0<T<T^*\). In particular, the time-continuity properties in \eqref{x63} and \eqref{x64} extend to \(t=T^*\). Consequently, there exist
	\[
	\rho^*\in H^3,
	\qquad
	u^*\in H^2,
	\qquad
	v^*\in H^2,
	\]
	such that
	\begin{equation}\label{b1}
	\rho(t)\to \rho^*
	\quad\text{in }H^3,
	\quad
	u(t)\to u^*
	\quad\text{in }H^2,
	\quad
	v(t)\to v^*
	\quad\text{in }H^2,
	\quad
	\text{as }t \to T^* .
	\end{equation}
	By \eqref{x62} and the embedding \(H^3(\Omega)\hookrightarrow
	C^1(\overline\Omega)\), we also have
	\begin{align}\label{x66}
		0<c\le \rho^*(x)\le C,
		\qquad x\in\overline\Omega .
	\end{align}
    Moreover, \eqref{b1} and the continuity of the trace operators imply that \((\rho^*,u^*)\) satisfies the conditions \eqref{x15} and \eqref{x113} required by the local existence theory Lemma~\ref{L41}.
	
	Taking \((\rho^*,u^*)\) as initial data at time \(T^*\), Lemma~\ref{L41} provides some \(\delta>0\) for which the solution can be extended to \([T^*,T^*+\delta]\). This contradicts the maximality of \(T^*\). Hence
	\[
	T^*=\infty .
	\]
	Consequently, the strong solution exists globally in time. Combining the estimates proved in Sections~3--5 with \eqref{x63}, \eqref{x64}, and \eqref{x65}, we obtain all the regularities and bounds stated in Theorem~\ref{th11} for \(N=3\).
	
	At this stage, it remains to deal with the 2D case. The proof is identical once one gets the upper bound of the density. The only difference is the range of \(\gamma\). When \(N=2\),
	the Sobolev embedding
	\[
	H^1(\Omega)\hookrightarrow L^p(\Omega),
	\qquad
	1\le p<\infty,
	\]
	together with the basic energy estimate yields
	\[
	\rho\in L^\infty(0,T;L^p),
	\qquad
	1\le p<\infty .
	\]
	Therefore, in the De Giorgi iteration for the upper bound of the density, the restriction on \(\gamma\) in 3D is no longer needed. Indeed, choosing \(q=2\), the relevant exponent satisfies
	\[
	q\gamma+2\gamma-q-1=4\gamma-3\ge1,
	\qquad
	\gamma\ge1.
	\]
	Thus, the upper bound of density holds for any \(\gamma\in[1,\infty)\). The remaining estimates then follow exactly the same steps as in the 3D case.
	Hence Theorem~\ref{th11} also holds for \(N=2\) and \(\gamma\in[1,\infty)\).
	
	For uniqueness, we refer to \cite{YuWu}. This completes the proof of Theorem~\ref{th11}.

    \section*{Acknowledgments}
    J. Li is supported by the National Natural Science Foundation of China (No. 12601413). B. L\"u is supported by  the National Natural Science Foundation of China (No.12371219),  and the Double-Thousand Plan of Jiangxi Province (No. jxsq2023201115).
    
    \section*{Conflict-of-interest statement}
    All authors declare that they have no conflicts of interest.

\end{document}